\documentclass[draft,11pt]{amsart}

\usepackage{amssymb}
\usepackage{amsmath, amsfonts,enumerate,amsxtra}
\usepackage{hyperref}
\usepackage{amsthm}
\usepackage{latexsym,bm}
\usepackage{euscript}

\let\cal\mathcal
\makeatletter \@mparswitchfalse \makeatother
\newtheorem{theorem}{Theorem}
\newtheorem{lemma}[theorem]{Lemma}

\newtheorem{corollary}[theorem]{Corollary}
\newtheorem{proposition}[theorem]{Proposition}
\theoremstyle{remark}
\newtheorem{remark}[theorem]{Remark}

\theoremstyle{definition}
\newtheorem{definition}[theorem]{Definition}
\newtheorem{problem}[theorem]{Problem}

\theoremstyle{remark}

\numberwithin{equation}{section}
\numberwithin{theorem}{section}

\def\Z{\mathbb Z}

\def\R{\mathbb R}
\def\M{\cal{M}}
\def\H{\cal{H}}
\let\<\langle
\def\ch{\raise 0.5ex \hbox{$\chi$}}
\def\T{\tau}
\def\E{\cal{E}}

\let\\\cr
\let\phi\varphi

\let\epsilon\varepsilon

\renewcommand{\a}{\alpha}
\renewcommand{\b}{\beta}
\newcommand{\g}{\gamma}

\renewcommand{\O}{{\Omega}}

\newcommand{\s}{\sigma}

\newcommand{\h}{\mathsf{h}}

\newcommand{\Tr}{\mbox{\rm Tr}}

\begin{document}

\title[Differentially subordinate martingales]{Square functions for noncommutative differentially subordinate  martingales}

\author[Jiao]{Yong Jiao}
\address{School of Mathematics  and Statistics, Central South University, Changsha 410075, China}
\email{jiaoyong@csu.edu.cn}

\author[Randrianantoanina]{Narcisse Randrianantoanina}
\address{Department of Mathematics, Miami University, Oxford,
Ohio 45056, USA}
 \email{randrin@miamioh.edu}
 
 \author[Wu]{Lian Wu}
\address{School of Mathematics  and Statistics, Central South University, Changsha 410075, China}
\email{wulian@csu.edu.cn}

\author[Zhou]{ Dejian Zhou}
\address{School of Mathematics  and Statistics, Central South University, Changsha 410075, China}
\email{zhoudejian@csu.edu.cn}

%\date{\today}
\thanks{Yong Jiao is supported by the NSFC  (No.11471337, No.11722114).  Lian Wu is supported by the NSFC (No.11601526).}

 
\begin{abstract} We  prove inequalities involving noncommutative differentially subordinate martingales. More precisely, we prove that if $x$ is a self-adjoint noncommutative martingale and $y$ is weakly differentially subordinate to $x$ then
$y$ admits a decomposition $dy=a +b +c$ (resp. $dy=z +w$) where $a$, $b$, and $c$ are adapted sequences (resp. $z$ and $w$ are martingale difference sequences) such that:
\[
\Big\| (a_n)_{n\geq 1}\Big\|_{L_{1,\infty}(\M\overline{\otimes}\ell_\infty)} +\Big\| \Big(\sum_{n\geq 1} \E_{n-1}|b_n|^2
\Big)^{{1}/{2}}\Big\|_{1, \infty} + \Big\|
\Big(\sum_{n\geq 1} \E_{n-1}|c_n^*|^2 \Big)^{{1}/{2}}\Big\|_{1,
\infty} \leq  C\big\| x \big\|_1
\]
(resp.
$$
\Big\| \Big(\sum_{n\geq1} |z_n|^2
\Big)^{{1}/{2}}\Big\|_{1, \infty} + \Big\|
\Big(\sum_{n\geq 1} |w_n^*|^2 \Big)^{{1}/{2}}\Big\|_{1,
\infty} \leq  C\big\| x \big\|_1).
$$
We also  prove strong-type $(p,p)$  versions of  the above weak-type  results for $1<p<2$.  In order to provide more insights into  
 the interactions  between  noncommutative  differential subordinations and martingale Hardy spaces when $1\leq p<2$, we also provide several martingale inequalities with sharp constants which are new and of independent interest.

As a byproduct of our approach, we obtain new and constructive proofs of both the noncommutative Burkholder-Gundy inequalities  and the noncommutative 
Burkholder/Rosenthal  inequalities for $1<p<2$ with the optimal order of the constants when $p \to 1$.

\end{abstract}

\maketitle

\section{Introduction}
It is a well known  fact that probabilistic inequalities and martingale inequalities in particular have broad impacts across many different fields of mathematics. Just like its commutative counterpart,  
 noncommutative martingale theory  has now  emerged  as a  very useful tool in various aspects of noncommutative analysis, noncommutative (or quantum) probability, and operator algebras. Recall that the origin of the current phase of development of  the theory of noncommutative martingale  comes from  the establishment  of the noncommutative Burkholder-Gundy inequality by Pisier and Xu in \cite{PX} . As explained in \cite{PX}, the interests on noncommutive martingales were primarily motivated by mathematical physics. 
 Since \cite{PX},  the theory of noncommutative martingale has been  steadily  progressing  to a point where many classical inequalities now have noncommutative analogues. The articles \cite{Bekjan-Chen-Perrin-Y, Hong-Junge-Parcet,Jiao-OS-Wu-2,Jiao-OS-Wu,Ju,Junge-Perrin,Junge-Zeng1,JX,OS1,PR,Perrin,Ran15}  contain samples of various noncommutative analogues of some of the most well known classical inequalities and techniques in the literature. We also refer to the book \cite[Chap.~14]{Pisier-Martingale} for  an up-to-date  overview of the  current status of the noncommutative martingale theory.   For the  classical  
 theory, 
  the so-called differential subordination occupies a prominent role.   The main objective of the present  paper is to further advance the topic of differential subordination in the noncommutative settings. To motivate our consideration, let us briefly describe the classical situation. Suppose that $(\O,\Sigma, \mathbb{P})$  is a probability space and  $f=(f_n)_{n\geq 1}$ and $g=(g_n)_{n\geq 1}$ are martingales. We say that $g$ is \emph{differentially subordinate} to $f$ if for every $n\geq 1$, the inequality $|dg_n| \leq |df_n|$ holds almost surely  (here,  $(df_n)_{n\geq 1}$ and $(dg_n)_{n\geq 1}$ are the martingale difference sequences of the martingales $f$ and $g$ respectively).  The notion of differential subordination  was introduced by Burkholder in \cite{Burkholder-diff} and became one of the fundamental tools  in martingale theory. To be more specific, it  
is  being used as general framework  for some basic operations in martingale theory such as martingale transforms and square functions. Two fundamental results proved by Burkhorder in \cite{Burkholder-diff}  assert that if $g$ is differentially subordinate to $f$ then we have the weak-type inequality 
\begin{equation}\label{equation1}
\|g\|_{1,\infty} \leq 2\|f\|_1
\end{equation}
and the strong type  $L_p$-bound
\begin{equation}\label{equation2}
\|g\|_p \leq (p^*-1)\|f\|_p, \quad 1<p<\infty,
\end{equation}
where  $p^*=\max\{p,p'\}$ with $p'$ being the conjugate index of $p$.
For more information on classical differential subordinations, we refer to the monograph \cite{OS}.

In  the recent article \cite{Jiao-OS-Wu}, the first and  third authors and Os\c ekowski 
thoroughly examined 
possible generalizations  of  the concept of differential subordination  in  the context of noncommutative martingales. As  it is often the case when dealing with the noncommutative case,  it turns out that one needs to work with two different versions of domination relations  according to $1\leq p<2$ or $p\geq 2$. One called  \emph{weak differential subordination} is needed in order to handle inequalities in the range $1\leq p<2$ while a weaker version called \emph{very weak differential subordination} is sufficient for the case $p\geq 2$. We refer to the preliminary  section below  for the exact formulations of these two notions of dominations.  The main achievement   in  the paper \cite{Jiao-OS-Wu} is twofold: the first is to identify the right formulations of noncommutative dominations  for the two separate cases described above and the second is to show  that  under  these appropriate dominations, the two inequalities \eqref{equation1} and \eqref{equation2} remain valid (but with different constants) for the noncommutative setting. Versions of differential subordinations were also considered in \cite{Jiao-OS-Wu-3} for noncommutative submartingales. Motivated by these results, we consider  in this paper the cases of square functions and conditioned square functions of differentially subordinate martingales. That is, 
estimating Hardy space norms of noncommutative differentially  subordinate martingales. We should emphasize that  the general theme considered here has trivial answer for the classical situation. Take for instance   the case of square functions:
 if $f=(f_n)_{n\geq 1}$ is a martingale on  a  given probability space  and $g=(g_n)_{n\geq 1}$ is differentially subordinate to $f$ then  $S(g)\leq S(f)$ where $S(f)$ and $S(g)$ refer to the  square functions of $f$ and $g$ respectively. By  the definition of classical martingale Hardy space $\H_p(\Omega)$, one immediately gets
 \begin{equation}\label{Hardy-1}
 \big\| g\big\|_{\H_p(\Omega)} \leq \big\| f\big\|_{\H_p(\Omega)}, \quad \text{for}\ 1\leq p\leq \infty.
 \end{equation}
By the  classical Burkholder-Gundy inequalities (\cite{BG}),  one readily gets that  for some constant $c_p$ 
\begin{equation}\label{Hardy-2}
\big\|g \big\|_{\H_p(\Omega)} \leq c_p \big\|f \big\|_p, \quad \text{for}\ 1<p<\infty.
\end{equation}
 Similarly,
 from a classical result of Burkholder on weak-type $(1,1)$ boundedness of square functions (\cite{Bu4}), one also easily deduces  that  
 \begin{equation}\label{weak-1}
 \|S(g)\|_{1,\infty} \leq \|S(f)\|_{1,\infty}\leq 2\|f\|_1.
 \end{equation} 
 Similar type  inequalities also hold for conditioned Hardy spaces norms and weak-type inequality involving conditioned square functions (see Section~4 below for more details).  

The noncommutative situation is radically different for the case $1\leq p<2$.
The main difference  lies with  the fact that in the noncommutative situation, we have two types of square functions  and  noncommutative martingale Hardy spaces 
consisting of sum of  row Hardy spaces and column Hardy spaces  when $1\leq p<2$.
As a result,   suitable decompositions are needed  when computing Hardy space norms. This phenomenon reveals that  noncommutative analogue of the weak-type inequality \eqref{weak-1} and  noncommutative analogues of \eqref{Hardy-1} and \eqref{Hardy-2} for $1\leq p<2$  become highly nontrivial. To further support this claim, assume that a noncommutative  martingale $y=(y_n)_{n\geq 1}$ is weakly differentially subordinate to another martingale $x=(x_n)_{n\geq 1}$.  Although we have  the noncommutative analogue of the Burkholder weak-type $(1,1)$ boundedness of square functions in \cite{Ran18}, it is  in the form of a decomposition $x=a +b$ such that $\|S_c(a)\|_{1,\infty} +
\|S_r(b)\|_{1,\infty} \leq C\|x\|_1$ where $S_c(\cdot)$ (resp.  $S_r(\cdot)$) denotes the
column (resp. row) square functions. Such decomposition is usually made up of non self-adjoint martingales and  it is very unclear if we can relate $y$ (or any decomposition of $y$) to $a$ and $b$ through some form of dominations. Therefore, a noncommutative analogue of \eqref{weak-1} cannot   be  easily  deduced. 
 It is our intent in this paper to clarify  this situation. In particular, we consider the question of whether  noncommutative analogues of \eqref{Hardy-1}, \eqref{Hardy-2}, \eqref{weak-1}, as well as their conditioned versions exist. As we will see below, only the case $1\leq p<2$   is of interest since when $p\geq 2$, the martingale Hardy space norms do not require any decomposition and therefore inequalities \eqref{Hardy-1} and \eqref{Hardy-2} are clearly satisfied by self-adjoint martingales under the assumption of very weak differential subordinations. Our approach for the weak-type situation was to consider concrete decompositions of weakly differentially subordinate martingales  in the spirit of  the decompositions used in  \cite{Ran18, Ran21}. 
More precisely,  we show (Theorem~\ref{main-weak} and Theorem~\ref{main-weak-S}) using concrete decompositions that noncommutative weak-type inequalities analogous to \eqref{weak-1} hold for mixture of column and  row conditioned square functions together with a  diagonal part  in the spirit of  \cite{Ran21} as well as mixture of column and row square functions formulated in the style of \cite{Ran18}. 
For strong-type $(p,p)$, we  establish the exact analogue of \eqref{Hardy-2} for $1<p<2$. 
  
  The paper is organized as follows.
In the next section, we  present some necessary background on noncommutative spaces and gather some  basic facts and preliminary results  concerning noncommutative martingales that we will need throughout. 

Section~3 includes  a new  description of a  Gundy type decomposition  for differentially subordinate martingales which is different from the versions in \cite{PR} and could be of independent interest. This new decomposition  is crucial in our approach to the weak-type $(1,1)$ situation. In fact,  it allows us to extend some of the techniques  used in \cite{PR}  for $L_1$-bounded martingales to weakly differentially subordinate martingales.
For instance, 
 using our version of Gundy's decomposition, we obtain a 
new proof of the weak-type $(1,1)$ for weakly  differentially subordinate  martingale \cite[Theorem~4.1]{Jiao-OS-Wu} in the same spirit as the proof of the  weak-type $(1,1)$  boundedness of noncommutative martingale transforms presented in \cite{PR}.

Section~4  contains  our principal results. More precisely, we present full descriptions of the two concrete decompositions of weakly differentially subordinate martingales and  show that  they satisfy weak-type (1,1) inequalities analogous to \eqref{weak-1}.
We also extrapolate that noncommutative analogues of \eqref{Hardy-2} hold when $1<p<2$ and the noncommutative Hardy spaces $\h_p(\M)$ and $\H_p(\M)$ are used. This section also contains noncommutative extensions of some sharp constant results due to Wang  in \cite{Wang}  
(Theorem~\ref{L-h}). To the best of our knowledge, sharp constant inequalities for noncommutative martingales have not been considered previously.
 Based on these noncommutative extensions of  Wang's results, our noncommutative analogue of \eqref{Hardy-2} using $\h_p(\M)$ (see Theorem~\ref{strong-type-1} below) implies the strong-type $(p,p)$ inequality from \cite[Theorem~5.1(i)]{Jiao-OS-Wu}. It turns out that our Theorem~\ref{L-h} and the previously described strong type $(p,p)$  results provide new and constructive proofs of both the noncommutative Burkholder-Gundy inequalities and the noncommutative Burkholder/Rosenthal inequalities when $1<p<2$.

In the last section, we discuss how  some  estimates from Section~4 can be used to  prove   some  moment inequalities associated with convex functions. The results  obtained in this section partially answer some open problems from \cite{Bekjan-Chen-2}.

%%%%%%%%%%%%
%%%%%%%%%%%%%

\section{Preliminaries}
\subsection{Noncommutative spaces}
Throughout this paper, $\M$ will always denote a finite  von Neumann algebra equipped with a normal  faithful normalized  trace $\T$. If $\M$ is acting on a Hilbert space $H$ then  a closed densely defined operator $x$  on $H$ is said to be affiliated with $\M$ if $u^*xu=x$ for all unitary operators $u$ in the commutant $\M'$ of $\M$.  If  $x$ is a densely defined self-adjoint  operator on $H$  and  $x=\int_{-\infty}^\infty sde_s^x$  is  its spectral decomposition then for any Borel subset $B \subseteq \R$, we denote by $\ch_B(x)$ the  corresponding spectral projection
$\int_{-\infty}^\infty \ch_B(s)\  de_s^x$.  Since $\M$ is finite, every
 closed densely defined operator $x$ affiliated with $\M$ is  $\T$-measurable in the sense that for every $\epsilon>0$, there exists a projection $p \in \M$ with $\T({\bf  1}-p)<\epsilon$ and $xp \in \M$.
  For a $\T$-measurable operator $x$, the decreasing function  on $[0,\infty)$ defined by $s\mapsto \T(\ch_{(s,\infty)}(|x|))$ will be referred to as the \emph{distribution function} of $x$. Denote by $L_0(\M,\T)$ the set of all $\T$-measurable operators. The set $L_0(\M,\T)$ is a $*$-algebra with respect to the strong sum, the strong product, and the adjoint operation. For $x \in L_0(\M,\T)$, the generalized singular-value function of $\mu(x)$ of $x$ is defined by:
\[
\mu_t(x)=\inf\big\{s\geq 0 : \T(\ch_{(s,\infty)}(|x|))\leq t\big\} \quad \text{for}\ t>0.
\] 
For a complete study of  generalized singular value functions and distributions functions,
we refer to \cite{FK}.  For the case where $\M$ is the abelian von Neumann algebra $L_\infty(0,1)$ with  the trace given by  integration  with respect to the Lebesgue  measure on $(0,1)$, $L_0(\M,\T)$  becomes  the linear space consisting of  those  measurable  complex functions  on $(0,1)$ which are bounded except on a set of arbitrarily small  measure and for  $f\in L_0(\M,\T)$,  $\mu(f)$ is the decreasing rearrangement of  the function $|f|$ in the sense of \cite{LT}.
For $0<p \leq \infty$, we denote by $L_p(\M,\T)$ or simply $L_p(\M)$ the noncommutative $L_p$-space associated with  the pair $(\M,\T)$. As usual, $L_\infty(\M,\T)$ is just the von Neumann algebra $\M$ with the operator norm. Beside the $L_p$-spaces,  we will also need to work with more general noncommutative symmetric spaces. A (quasi) Banach function space  $(E,\|\cdot\|_E)$ of measurable functions  on the interval $(0,1)$ is called \emph{symmetric} if for any $g \in E$ and any $f \in L_0(0,1))$ with $\mu(f) \leq \mu(g)$, we have $f \in E$ and $\|f\|_E \leq \|g\|_E$.

For a given symmetric (quasi) Banach function space  $(E, \|\cdot\|_E)$ on  $(0, 1)$, we define the  corresponding  noncommutative space by setting:
\begin{equation*}
E(\M, \T) = \big\{ x \in
L_0(\M,\T)\ : \ \mu(x) \in E \big\}. 
\end{equation*}
Equipped with the (quasi)  norm
$\|x\|_{E(\M,\T)} := \| \mu(x)\|_E$, the space  $E(\M,\T)$ is a complex (quasi)  Banach space (\cite{Kalton-Sukochev,X}) and is usually referred to as the \emph{noncommutative symmetric space} associated with $(\M,\T)$ corresponding to  the symmetric  space $(E, \|\cdot\|_E)$.  When $E=L_p(0,1)$ for some $0<p<\infty$, then $E(\M,\T)$ coincides with $L_p(\M,\T)$.
The particular case of noncommutative  weak-$L_1$ space  $L_{1,\infty}(\M,\T)$ will be heavily used. This is defined as the collection of all  $x \in L_0(\M,\T)$ for which the quasi-norm
\[
\big\|x \big\|_{1,\infty} = \sup_{t>0} t\mu_t(x)=\sup_{\lambda >0} \lambda \T\big( \ch_{(\lambda,\infty)}(|x|)\big)
\]  
is finite. According to the general construction described above, the linear space $L_{1,\infty}(\M,\T)$ is a quasi-Banach space. In the sequel, we will also use weak-$L_1$ space associated with  the semifinite von neuamann algebra $\M\overline{\otimes} \ell_\infty$  whose definition is identical to the finite case.
 We refer to the survey \cite{PX3} for more in depth  treatment of noncommutative spaces.
We   end this subsection by recording  a  general
 inequality  on distribution  function of operators    that we will  use repeatedly in the sequel.  
\begin{lemma} \label{Quasi-Triangle}
Let  $a$ and $b$ be $\T$-measurable operators and  $t, s > 0$. We have 
\[
 \T \Big(
\ch_{(t+s,\infty)} \big( |a+b| \big) \Big) \leq  
\T \Big( \ch_{(t, \infty)} \big( |a| \big) \Big) +
  \T \Big( \ch_{(s, \infty)} \big( |b|
\big) \Big).
\]
\end{lemma}
\begin{proof}
We note first that $|a+b| \leq u|a|u^* +v|b|v^*$ where  $u$ and $v$ are  partial isometries from $\M$. Since  $\ch_{(t+s, \infty)}(|a+b|)$ is equivalent to a subprojection of 
$\ch_{(t+s, \infty)}(u|a|u^* +v|b|v^*)$, we have 
\[
\T(\ch_{(t+s, \infty)}(|a+b|)) 
\leq \T(\ch_{(t+s, \infty)}(u|a|u^* +v|b|v^*)).
\]
 Next, we note that since  $u|a|u^*$ and $v|b|v^*$ are positive operators, we further get  according to 
 \cite[Lemma~16]{Sukochev-Zanin} that
\[
\T(\ch_{(t+s, \infty)}(|a+b|)) 
\leq \T(\ch_{(t, \infty)}(u|a|u^*)) + \T(\ch_{(s, \infty)}(v|b|v^*)).
\] 
To conclude the proof, we make the simple observation that 
\[\T \big( \ch_{(t, \infty)} ( u|a|u^* ) \big)\leq  \T \big( \ch_{(t, \infty)} ( |a| ) \big)  \ \text{and} \ 
 \T \big( \ch_{(s, \infty)} (v |b|v^*
) \big) \leq  \T \big( \ch_{(s, \infty)} (|b|)
 \big)\]
 These  follow easily from the property of distribution functions  that for any given $\T$-measurable operator $x$,  the identity  $\T \big( \ch_{(t, \infty)} ( x^*x ) \big)=\T \big( \ch_{(t, \infty)} ( xx^* ) \big)$ holds. Indeed, using $x=|a|^{1/2} u^*$, we immediately obtain that $x^*x=u|a|u^*$ and $xx^*\leq |a|$, thus the first inequality follows. Identical  arguments can be applied to $v|b|v^*$.
\end{proof}
\subsection{Noncommutative martingales}

In this subsection,  we will review the basics of noncommutative martingales, recall some recently introduced notions of noncommutative differential subordinations for martingales,  and present some preliminary results that we will need in the sequel. 

\subsubsection{Definitions and martingale Hardy spaces}
Let
$(\M_n)_{n\geq1}$ be an increasing sequence of von Neumann
subalgebras of $\M$ such that the union of the $\M_n$'s is
$w^*$-dense in $\M$. Since $\M$ is finite,  for every $n\geq 1$, there exists a $\T$-invariant   conditional
expectation from  $\M$ onto  $\M_n$.   Since $\E_n$ is $\T$-invariant, it extends to a contractive projection from $L_p(\M,\T)$ onto $L_p(\M_n,\T_n)$ for all $1\leq p\leq \infty$, where $\T_n$ denotes the restriction of $\T$ on $\M_n$.

A sequence $x=(x_n)_{n\geq 1}$ in
$L_1(\M)$ is called  \emph{a noncommutative martingale} with respect to the filtration 
$(\M_n)_{n\geq1}$ if  for every $n\geq 1$,
\[
\E_n(x_{n+1})=x_n.\]
 If in addition, all the $x_n$'s are in $L_p(\M)$ for some
$1\leq p\leq\infty$, $x$ is called an $L_p$-martingale. In this case, 
we set:
\[
\|x\|_p=\sup_{n\geq1}\|x_n\|_p.\]
If $\|x\|_p<\infty$, $x$ is called an  $L_p$-bounded martingale.
For $n\geq 1$, we  define $dx_n=x_n-x_{n-1}$  with the
convention that $x_0=0$ and $\mathcal{E}_0=\mathcal{E}_1$. The
sequence $dx=(dx_n)_{n\geq 1}$ is called the martingale difference sequence
of  the martingale $x$.  The martingale $(x_n)_{n\geq 1}$ is said to be finite if there exists $N \geq 1$ such that $x_n=x_N$ for all $n\geq N$.

We now  review  the construction of various Hardy spaces for noncommutative martingales. We begin with  descriptions to some   general  spaces   that we will need in the sequel. For $0<p \leq \infty$ 
and a  sequence $a=(a_n)_{n\geq 1}$ in $L_p(\M)$, we set 
\[
\big\|a\big\|_{L_p(\M;\ell_2^c)}=\Big\| \big(\sum_{n\geq 1} |a_n|^2\big)^{1/2}\Big\|_p.
\] 
We  define $L_p(\M;\ell_2^c)$ to be the collection of  all sequence $a=(a_n)_{n\geq 1}$ in $L_p(\M)$ for which the quantity  $\big\|a\big\|_{L_p(\M;\ell_2^c)}$ is finite. It is well-known that  when  equipped with $\|\cdot\|_{L_p(\M;\ell_2^c)}$, the linear space $L_p(\M;\ell_2^c)$ becomes a (quasi)-Banach space. We refer to \cite{PX, PX3} for this fact.
 We   will also need  the conditioned $L_p$-spaces which are defined as follows:
for  $0< p \leq \infty$ and  a   sequence $a=(a_n)_{n\geq 1}$ in $ L_2(\M)$, we set
\[
\big\|a\big\|_{L_p^{\rm cond}(\M;\ell_2^c)} = \Big\| \big( \sum_{n\geq 1} \E_{n-1} (a_n^* a_n) \big)^{1/2}\Big\|_p.
\]
For $0<p<2$,  we define $L_p^{\rm cond}(\M;\ell_2^c)$ to be the completion of the space of all finite  sequences in $ L_2(\M)$ equipped with the (quasi) norm $\| \cdot\|_{L_p^{\rm cond}(\M;\ell_2^c)} $ while for $2\leq  p \leq \infty$, we may define $L_p^{\rm cond}(\M;\ell_2^c)$ directly as the set of all sequences $a=(a_n)_{n\geq1}$  in $L_p(\M)$ for which the increasing sequence $\big\{\big(\sum_{k=1}^n \E_{k-1}(a_k^*a_k) \big)^{1/2}\big\}_{n\geq 1}$ is bounded in $L_p(\M)$. In this range,  $L_p^{\rm cond}(\M;\ell_2^c)$ is also equipped with the norm $\| \cdot\|_{L_p^{\rm cond}(\M;\ell_2^c)} $. We can  extend the above definition in the context of weak-$L_1$-space by setting for any sequence $a=(a_n)_{n\geq 1}$ in $ L_2(\M)$, 
\[
\big\|a\big\|_{L_{1,\infty}^{\rm cond}(\M;\ell_2^c)} = \Big\| \big( \sum_{n\geq 1} \E_{n-1} (a_n^* a_n) \big)^{1/2}\Big\|_{1,\infty}.
\]
We refer to \cite{Ju,JX} for more details on these conditioned spaces.
Following \cite{PX}, we define  the  column  and row versions of square functions
relative to a  martingale $x = (x_n)_{n\geq 1}$:
 \[
 S_{c,n} (x) = \Big ( \sum^n_{k = 1} |dx_k |^2 \Big )^{1/2}, \quad
 S_c (x) = \Big ( \sum^{\infty}_{k = 1} |dx_k |^2 \Big )^{1/2}\,.
 \]
 and 
 \[
 S_{r,n} (x) = \Big ( \sum^n_{k = 1} |dx_k^* |^2 \Big )^{1/2}, \quad
 S_r (x) = \Big ( \sum^{\infty}_{k = 1} |dx_k^* |^2 \Big )^{1/2}\,.
 \]
For $1\leq p \leq \infty$, the column  martingale Hardy space 
 $\H_p^c (\M)$ (resp. the row martingale Hardy space $\H_p^r(\M)$) is defined to be the space of all martingales $x$ for which $S_c(x) \in L_p(\M)$ (resp. $S_r(x) \in L_p(\M)$) under the norm, $\|x\|_{\H_p^c}=\|S_c(x)\|_p$ (resp. $\|x\|_{\H_p^r}=\|S_r(x)\|_p$).  For $0<p<1$, $\H_p^c(\M)$ (resp. $\H_p^r(\M)$) is the completion of all finite martingale $x \in L_2(\M)$ under the quasi-norm $\|\cdot\|_{\H_p^c}$ (resp. $\|\cdot\|_{\H_p^r}$).
 The noncommutative martingale Hardy spaces $\H_p(\M)$  are
defined as follows: if $0< p<2$,
\[\H_p(\M)=\H^c_p(\M)+\H^r_p(\M)
\]
equipped with the (quasi) norm
\[
\|x\|_{\H_p}=\inf\big\{\|y\|_{\H^c_p}+\|z\|_{\H^r_p}\big\}
\]
where the infimum is taken over all decomposition $x=y +z$ with $y\in\H_p^c$ and $z \in \H_p^r(\M)$.
When $2\leq p<\infty$,
$$\mathcal{H}_p(\M)=\mathcal{H}^c_p(\M)\cap\mathcal{H}^r_p(\M)$$
equipped with the norm
$$\|x\|_{\mathcal{H}_p}=\max\big\{\|x\|_{\mathcal{H}^c_p},\|x\|_{\mathcal{H}^r_p}\big\}.$$

We now consider the conditioned version of $\H_p^c(\M)$  and $\H_p^r(\M)$ developed in \cite{JX}.
Let $x=(x_n)_{n \geq 1}$ be a martingale in $L_2(\M)$.
We set
 \[
 s_{c,n} (x) = \Big ( \sum^n_{k = 1} \E_{k-1}|dx_k |^2 \Big )^{1/2}, \quad
 s_c (x) = \Big ( \sum^{\infty}_{k = 1} \E_{k-1}|dx_k |^2 \Big )^{1/2}\,.
 \]
 The operator $s_c(x)$ is  called  the column conditioned square function of $x$.
 For convenience, we will use the notation 
 \[
  \s_{c,n} (a) = \Big ( \sum^{n}_{k = 1} \E_{k-1}|a_k |^2 \Big )^{1/2}, \quad
 \s_c (a) = \Big ( \sum^{\infty}_{k = 1} \E_{k-1}|a_k |^2 \Big )^{1/2}
 \]
 for sequence $a=(a_k)_{k\geq 1}$ in $L_2(\M)$  that is not necessarily a martingale difference sequence.
For $2\leq p \leq \infty$, the column conditioned martingale Hardy space 
$\h_p^c(\M)$  is defined to be  the space of all martingales $x$ for which $s_c(x)$ belongs to $L_p(\M)$, equipped with the norm  $\| x \|_{\h_p^c}=\| s_c (x) \|_p$.
We refer to \cite{JX} for the fact $(\h_p^c(\M), \|\cdot\|_{\h_p^c})$  is a Banach space. 
For $0<p<2$, we
define $\h_p^c (\M)$ to be  the completion of the linear space of finite martingales in $L_2(\M)$ under the (quasi) norm $\| x \|_{\h_p^c}=\| s_c (x) \|_p$. Obvious modification as before  is made to describe the row versions.
We will also need a third type of Hardy space known as the diagonal Hardy space $\h^d_p(\M)$. This is defined 
as the subspace of $\ell_p(L_p(\M))$ consisting of all martingale
difference sequences.  The conditioned  Hardy spaces are defined
as follows: if $0<p<2$,
\[
\h_p(\M)=\h^c_p(\M)+\h^r_p(\M)+\h^d_p(\M)\]
equipped with the (quasi) norm
\[
\|x\|_{\h_p}=\inf\big\{\|y\|_{\h^c_p}+\|z\|_{\h^r_p}+\|w\|_{\h^d_p}\big\}\]
where the infimum is taken over all decomposition $x= y+z+w$ with $y \in \h_p^c(\M)$, $z \in \h_p^r(\M)$, and $z \in \h_p^d(\M)$.  When $2\leq p<\infty$,
\[
\mathsf{h}_p(\M)=\h^c_p(\M)\cap\h^r_p(\M)\cap\h^d_p(\M)
\]
equipped with the norm
\[
\|x\|_{\h_p}=\max\big\{\|x\|_{\h^c_p},\|x\|_{\h^r_p},\|x\|_{\h^d_p}\big\}.\]

\subsubsection{Differential subordination}

We now isolate three different types of  differential subordination of  noncommutative martingales which constitute  the main focus of the present paper. 

\begin{definition}[\cite{OS1}]\label{subordinate1}
Let $x$, $y$ be two self-adjoint $L_2$-martingales. We say that $y$ is \emph{ differentially subordinate} to $x$ if the following two conditions hold:
\begin{enumerate}[{\rm(i)}]
\item for any $n\geq 1$ and any projection $R \in \M_n$, we have
\begin{equation*}
\T(Rdy_nRdy_n R) \leq  \T(Rdx_nRdx_n R);
\end{equation*}
\item for any $n\geq 1$ and any orthogonal projections $R$, $S$ in $\M_n$ such that $R+S \in \M_{n-1}$, we have
\begin{equation*}
\T(Rdy_nSdy_n R) \leq  \T(Rdx_nSdx_n R);
\end{equation*}
\end{enumerate}
\end{definition}

\begin{definition}[\cite{Jiao-OS-Wu}]\label{subordinate2}
Let $x$, $y$ be two self-adjoint martingales. We say that $y$ is \emph{weakly differentially subordinate} to $x$ if for any $n\geq 1$ and any projection $R \in \M_{n-1}$, we have
\begin{equation}\label{w-sub}
Rdy_nRdy_n R \leq  Rdx_nRdx_n R.
\end{equation}
We say that $y$ is \emph{very weakly differentially subordinate} to $x$ if for every $n\geq 1$, we have 
\begin{equation}\label{vw-sub}
dy_n^2 \leq dx_n^2.
\end{equation}
\end{definition}
In the commutative case, all three dominations are equivalent to  the property that 
$|dy_n| \leq |dx_n|$ for all $n\geq 1$, which is the classical notion of differential subordination introduced by Burkholder in \cite{Burkholder-diff}.  In the noncommutative  setting, it is clear from the definitions that the differential subordination implies the the weak differential subordination  and the latter  implies the very weak differential subordination. However,  it was shown in \cite{Jiao-OS-Wu} that the three notions are not equivalent in general.  The main discovery of  \cite{Jiao-OS-Wu} is that the  notion of weak differential  subordination is   needed  in order to generalize weak-type $(1,1)$ and strong-type $(p,p)$ results from the classical setting to the noncommutative setting  when $1<p<2$, while the notion of very weak differential subordination is sufficient for the range $p\geq 2$. 
Below, we will concentrate on the range $1\leq p<2$ and  investigate how weak  differential subordinations interact with noncommutative martingale Hardy spaces.  
We note  that martingale transforms with commuting symbols (in the sense of \cite{PX,Ran15}) of self-adjoint martingales are examples of weak differential subordinations. Indeed, if $(x_n)_{n\ge 1}$ is a self-adjoint martingale and for each $n\geq 2$, $\xi_{n-1}$ is a  self-adjoint contraction that belongs to $\M_{n-1} \cap \M_n'$, then by commutation, one can easily verify that
for every  projection $R \in \M_{n-1}$,
\[
R\xi_{n-1}dx_n R\xi_{n-1}dx_n R \leq Rdx_nRdx_nR.
\]
For the case where the  $\xi_{n-1}$'s are not self-adjoint, one can consider their real and imaginary parts  separately and get combinations of two weakly differentially subordinate martingales.
Another example of weak differential subordinate martingale was also exhibited in the proof of \cite[Lemma~3.3]{Jiao-OS-Wu}.

\subsubsection{Cuculescu projections and their relatives}

We recall  the so-called Cuculescu projections associated to a given  self-adjoint $L_1$-bounded martingale $x$. 
Such sequence of projections will
 play crucial role in the construction below. We will also introduce
 some other sequences of projections derived from the Cuculescu projections and gather some of their  properties that are relevant  for our proofs.

  Fix $\lambda$ to be a positive real number. Set $q_0^{(\lambda)}={\bf 1}$ and inductively we define the decreasing sequence of projections
\[
q_n^{(\lambda)} := q_{n-1}^{(\lambda)} \ch_{[-\lambda,\lambda]} \Big(
q_{n-1}^{(\lambda)} x_{n} q_{n-1}^{(\lambda)} \Big)=\ch_{[-\lambda,\lambda]}\Big(
q_{n-1}^{(\lambda)} x_{n} q_{n-1}^{(\lambda)} \Big)q_{n-1}^{(\lambda)} .
\]
The sequence $(q_n^{(\lambda)})_{n\geq 1}$ was first considered in \cite{CUC} and will be referred to as the Cuculescu projections associated with $\lambda$. Their significance  in  the area of noncommutative martingales  is now well-established. We record  some of the basic properties that we will use.
\begin{proposition}[{\cite[Proposition~1.4]{PR}}]\label{Cuculescu}
The sequence $(q_n^{(\lambda)})$ satisfies the following properties:
\begin{itemize}
\item[(i)] for every $n \geq 1$, $q_{n}^{(\lambda)} \in \M_{n}$;
\item[(ii)] for every $n \geq 1$, $q_{n}^{(\lambda)}$ commutes
with $q_{n-1}^{(\lambda)} x_{n} q_{n-1}^{(\lambda)}$;
\item[(iii)] for every $n \geq 1$, $|q_{n}^{(\lambda)} x_{n}
q_{n}^{(\lambda)}| \leq \lambda q_{n}^{(\lambda)}$. In particular, $\|q_{n}^{(\lambda)} x_{n}
q_{n}^{(\lambda)}\|_\infty \leq \lambda$;
\item[(iv)]
for every $N\geq 1$, 
\[
\T \Big( {\bf 1}-q_N^{(\lambda)} \Big)
\leq \frac{1}{\lambda} \T \Big( ({\bf 1}-q_N^{(\lambda)}) |x_N|\Big)\leq \frac{1}{\lambda} \|x\|_1. 
\]
\end{itemize}
\end{proposition}

  Following
\cite{Ran18,Ran21}, we consider   collection of projections derived from the Cuculescu projections. For $n\geq 1$ and $i \in \Z$, we set:
\[
e_{i,n} := \bigwedge^\infty_{k=i}q_{n}^{(2^k)} \ \ \text{and}\ \ \pi_{i,n}:=e_{i,n} -e_{i-1,n}.
\]
The family $\{e_{i,n}\}_{n\geq 1, i\in \Z}$ is decreasing on $n$ and increasing on $i$. Therefore,  for every $n\geq 1$, $(\pi_{i,n})_{i\in \Z}$ is a sequence of pairwise disjoint projections satisfying the  trivial but crucial identity  that  for every $k\in \Z$,
\[
\sum_{i=-\infty}^k  \pi_{i,n}=e_{k,n}
\]
where the convergence of the series is relative to the strong operator topology.
We also note that since ${\bf 1}-e_{k,n} \to 0$ when $k \to \infty$, we have 
$\sum_{i\in \Z}  \pi_{i,n}={\bf 1}$ for the strong operator topology.

 The family  of projections $(e_{k,n})_{k,n}$ satisfies similar properties as displayed by $(q_n^{(\lambda)})$ in Proposition~\ref{Cuculescu}(iv). More precisely, for $N\geq 1$ and $k\in \Z$, the following inequality holds:
\begin{equation}\label{e}
\T\big({\bf 1}-e_{k,N}\big) \leq \frac{1}{2^{k-1}} \T\big( ({\bf 1}-e_{k,N})|x_N| \big).
\end{equation}
To see this, we have by  the definition  of $e_{k,N}$ and Proposition~\ref{Cuculescu} that 
\[
\T\big({\bf 1}-e_{k,N}\big) \leq \sum_{j\geq k} \T\big({\bf 1}-q_N^{(2^j)}\big) \leq \sum_{j\geq k} 2^{-j}\T\big(({\bf 1}-q_N^{(2^j)} )|x_N|\big).
\]
Since for $j\geq k$, $q_N^{(2^j)} \geq e_{j,N} \geq e_{k,N}$, we have 
\begin{align*}
\T\big({\bf 1}-e_{k,N}\big)  &\leq (\sum_{j\geq k} 2^{-j})\T\big(({\bf 1}-e_{k,N})|x_N|\big)\\
&=2^{-k+1}\T\big(({\bf 1}-e_{k,N})|x_N|\big).
\end{align*}

In the sequel, we will  mainly use the corresponding sequence of pairwise disjoint projections  by grouping together the $\pi_{i,n}$'s when $i\leq 0$. That is, 
for $n\geq 1$, we set:
\begin{equation*}
\begin{cases}
p_{0,n} &:=\displaystyle{ e_{0,n}} \\
p_{i,n} &:=\displaystyle{ \pi_{i,n} \quad \text{for $i\geq
1$}}.
\end{cases}
\end{equation*}
 Then,
we have the following  basic properties:
\begin{enumerate}[$\bullet$]
\item  For any  given $n\geq 1$, $(p_{k,n})_{k\geq 0}$  is  a sequence of pairwise disjoint projections in $\M_{n}$.
\item For any  given $m\geq 1$,  $\sum_{k=0}^m  p_{k,n}=e_{m,n}$ and  $\sum_{k=0}^\infty p_{k,n}= {\bf 1}$ for the strong operator topology.
\end{enumerate}
These various family of projections  play important role in our construction in the next section. We now gather some auxiliary  inequalities that are essential in  our presentation.   
%%%%%

%%%%%%%%%%%

%%%%%%%%%%%%%%

\begin{lemma}\label{L2-norm} \begin{enumerate}[{\rm (i)}] 
\item For every $\lambda>0$, the following inequality holds:
\[
\sum_{n=1}^N \T\Big( q_n^{(\lambda)}dx_n q_{n-1}^{(\lambda)}dx_n q_n^{(\lambda)}\Big)\leq  \big\| q_N^{(\lambda )}x_Nq_N^{(\lambda )}\big\|_2^2 + 2\lambda \T\big( ({\bf 1}-q_{N}^{(\lambda)})|x_N|\big). 
\]

\item For any given  $k\geq 0$, the following inequality holds:
\[
\sum_{n= 1}^N \T\Big( e_{k,n}dx_n e_{k,n-1}dx_n e_{k,n}\Big)\leq 2^{k+1} \|x_N\|_1.
\]
\item For every $N\geq 1$,
\[
\sum_{n=1}^N \T\Big( e_{k,n}dx_n e_{k,n-1}dx_n e_{k,n}\Big) \leq  \big\| e_{k,N} x_N  e_{k,N}\big\|_2^2  + 6\ . \ 2^k \T\big( ({\bf 1}-e_{k,N} )|x_N|\big) .
\]
\end{enumerate}
\end{lemma}

\begin{proof} The first inequality is from \cite[Lemma~4.3]{Jiao-OS-Wu}. Below, we write for $m\geq 1$ and $k\geq 0$, $q_{k,m}$ for $q_{m}^{(2^k)}$.

For the second inequality, we make the observation that since $e_{k,m} \leq q_{k,m}$ for $m\geq 1$ and $k\geq 0$, we have
\[
\T\Big( e_{k,n}dx_n e_{k,n-1}dx_n e_{k,n}\Big) \leq  \T\Big( q_{k,n}dx_nq_{k,n-1} dx_n q_{k,n}\Big).
\]
It follows from the first inequality that:
\begin{equation}\label{first-estimate}
\sum_{n=1}^N  \T\Big( e_{k,n}dx_n e_{k,n-1}dx_n e_{k,n}\Big) \leq
\| q_{k,N}x_N  q_{k,N} \|_2^2 + 2^{k+1} \T\big( ({\bf 1}-q_{k,N})|x_N| \big).
\end{equation}
To deduce the second inequality, it suffices to observe that $\| q_{k,N}x_N  q_{k,N}\|_2^2  \leq 2^k \T(q_{k,N}|x_N|)$. Indeed,  $\| q_{k,N}x_N  q_{k,N}\|_2^2 \leq \|q_{k,N}x_N  q_{k,N}\|_\infty \|q_{k,N}x_N  q_{k,N}\|_1 \leq 2^k \|q_{k,N}x_N  q_{k,N}\|_1$ and by writing $x_N=x_N^+ -x_N^-$, it follows  from triangle inequality that $\|q_{k,N}x_N  q_{k,N}\|_1\leq \T(q_{k,N}|x_N|)$. We can then deduce from \eqref{first-estimate} that
\[
\sum_{n=1}^N  \T\Big( e_{k,n}dx_n e_{k,n-1}dx_n e_{k,n}\Big) \leq  2^k\big[ 2\T(|x_N|) -\T(q_{k,N}|x_N|) \big]\leq 2^{k+1} \T(|x_N|).
\]

\smallskip

In order to verify the third inequality,  we need  to majorize    $\| q_{k,N}x_N  q_{k,N} \|_2^2$ in \eqref{first-estimate}  in terms of  $\| e_{k,N}x_N  e_{k,N} \|_2^2$. Since $q_{k,N}-e_{k,N}$ and $e_{k,N}$ are two disjoint projections, we have
\begin{align*}
\| q_{k,N}x_N  q_{k,N} \|_2^2 &=\| (q_{k,N}-e_{k,N})q_{k,N}x_N  q_{k,N} \|_2^2 + \| e_{k,N}x_N q_{k,N} (q_{k,N}-e_{k,N})\|_2^2 +\| e_{k,N}x_N  e_{k,N} \|_2^2\\
&\leq  2^{2k+1} \T \big(q_{k,N} -e_{k,N}\big) + \| e_{k,N}x_N  e_{k,N} \|_2^2\\
&\leq  2^{2k+1}   \T \big({\bf 1} -e_{k,N}\big) +\| e_{k,N}x_N  e_{k,N} \|_2^2\\
&\leq  2^{k+2} \T\big( ({\bf 1} -e_{k,N})|x_N|\big) +\| e_{k,N}x_N  e_{k,N} \|_2^2
\end{align*}
where in the last inequality we use \eqref{e}. We obtain the inequality as stated by combining  this last  estimate with \eqref{first-estimate}. 
\end{proof}

The next two lemmas deal with martingales that are $L_p$-bounded  for some $p>1$.

\begin{lemma}\label{lem-p}  Let $1<p <\infty$ and $x=(x_n)_{n\geq 1}$ be a self-adjoint $L_p$-bounded martingale then for every $N\geq 1$,
\[
\sum_{j=0}^\infty 2^{(p-1)j} \T\big( ({\bf 1}- e_{j,N})|x_N|\big) \leq  \frac{2^{(p-1)^2}}{(2^{p-1}-1)^p}
  \big\|x_N\big\|_p^p.
\]
\end{lemma}
\begin{proof}
Since ${\bf 1}-e_{j, N} = \sum_{m \geq j+1}p_{m,N}$, we have
 \begin{align*}
 \sum_j 2^{(p-1)j} \T\big(({\bf 1} - e_{j,N})|x_N|\big)
 &= \sum_j 2^{(p-1)j} \sum_{m\geq j+1}\T\big(p_{m,N}|x_N| \big)\\
 &=\sum_{m\geq 1} \big(\sum_{j\leq m-1} 2^{(p-1)j} \big)\T\big(p_{m,N}|x_N| \big)\\
 &\leq \frac{1}{2^{p-1}-1}\, \T\Big( \big[\sum_{m\geq 1} 2^{(p-1)m} p_{m,N}]|x_N| \Big).
 \end{align*}
 We note that  if $x_N \in L_p(\M)$, then according to \cite[Lemma~5.3]{Jiao-OS-Wu},  the operator $\sum_{m\geq 1} 2^{(p-1)m} p_{m,N}$ belongs to $L_{p/(p-1)}(\M)$.  
 Using H\"older's inequality together with \eqref{e}, we get
 \begin{align*}
 \sum_j 2^{(p-1)j} \T\big(({\bf 1} - e_{j,N})|x_N|\big) &\leq \frac{1}{2^{p-1}-1}
 \Big(\sum_{m\geq 1} 2^{pm} \T\big(p_{m,N}\big)\Big)^{(p-1)/p} \big\|x_N\big\|_p\\
 &\leq\frac{1}{2^{p-1}-1}
 \Big(\frac{2^{p-1}}{2^{p-1} -1} \big\|x_N\big\|_p         \Big)^{p-1} \big\|x_N\big\|_p\\
 &=  \frac{2^{(p-1)^2}}{(2^{p-1}-1)^p}
  \big\|x_N\big\|_p^p
 \end{align*}
where in the second inequality, we use the estimate  from \cite[Lemma~5.3]{Jiao-OS-Wu}.
\end{proof}

\begin{lemma}\label{last} Let $1<p<2$ and assume that $x=(x_n)_{n\geq 1}$ is a self-adjoint  martingale that is $L_p$-bounded. Then for every  $N\geq 1$,
\[
\sum_{k=0}^\infty  2^{(p-2)k} \big\| e_{k,N} x_N e_{k,N}\big\|_2^2 \leq   \frac{2^{p^2 +1}}{(1-2^{p-2})(2^{p-1}-1)^p} \big\|x_N\big\|_p^p.
\]
\end{lemma}
\begin{proof} 
For $k\geq 0$, we claim that the following inequality holds:
\begin{equation}\label{est-1}
\big\|e_{k,N} x_N e_{k,N} \big\|_2^2 \leq 2 \sum_{j=-\infty}^k 2^{2j} \T\big( e_{j,N}-e_{j-1,N}\big).
\end{equation}
This is implicit in \cite{Jiao-OS-Wu} but we include the argument for completeness.
First, we recall that $e_{k,N}=\sum_{i=-\infty}^{k}\pi_{i,N}$ and the $\pi_{i,N}$'s  are pairwise  disjoint.  This implies that
\begin{align*}
\big\| e_{k,N} x_N e_{k,N}\big\|_2^2&=\sum_{-\infty< i,j\leq k}\|\pi_{i,N}x_N\pi_{j,N}\|_2^2\\
&\leq 2\sum_{-\infty< i\leq j\leq k}\|\pi_{i,N}x_N \pi_{j,N}\|_2^2\\
&=2\sum_{-\infty< j\leq k}\|e_{j,N}x_N \pi_{j,N}\|_2^2\\
&=2 \sum_{-\infty\leq j\leq k} \T(\pi_{j,N} [ e_{j,N} x_N e_{j,N} x_N e_{j,N}] \pi_{j,N} ).
\end{align*}
Since $\|e_{j,N} x_N e_{j,N}\|_\infty \leq 2^j$,  it follows that 
$
\big\| e_{k,N} x_N e_{k,N}\big\|_2^2 \leq  2\sum_{j=-\infty}^k 2^{2j}\T(\pi_{j,N})$. Thus, \eqref{est-1} is verified.

With this estimate, we may deduce that 
\begin{align*}
\sum_{k=0}^\infty2^{(p-2)k} \big\| e_{k,N} x_N e_{k,N}\big\|_2^2 &\leq  2\sum_{k=0}^\infty 2^{(p-2)k}\sum_{j=-\infty}^k 2^{2j} \T\big( e_{j,N}-e_{j-1,N}\big)\\
&=2\sum_{j=-\infty}^\infty 2^{2j}\Big(\sum_{k\geq j\vee 0} 2^{(p-2)k}\Big) \T\big( e_{j,N}-e_{j-1,N}\big)\\
&=\frac{2}{1-2^{p-2}}\Big[\sum_{j=-\infty}^0 2^{2j} \T\big( e_{j,N}-e_{j-1,N}\big) + \sum_{j=1}^\infty 2^{pj} \T\big( e_{j,N}-e_{j-1,N}\big)\Big]\\
&\leq \frac{2}{1-2^{p-2}}\sum_{j=-\infty}^\infty 2^{pj} \T\big( e_{j,N}-e_{j-1,N}\big)\\
&=\frac{2^{p+1}}{1-2^{p-2}} \sum_{j=-\infty}^\infty 2^{pj} \T\big( e_{j+1,N}-e_{j,N}\big).
\end{align*}
Let $a_N =  \sum_{j=-\infty}^\infty 2^{j} ( e_{j+1,N}-e_{j,N})$. By \cite[Lemma~5.3]{Jiao-OS-Wu}, $a_N \in L_p(\M)$ and satisfies:
\[
\|a_N\|_p \leq \frac{2^{p-1}}{2^{p-1}-1} \|x_N\|_p.
\]
Therefore, we arrive at the following estimate:
\begin{align*}
\sum_{k=0}^\infty2^{(p-2)k} \big\| e_{k,N} x_N e_{k,N}\big\|_2^2 &\leq \frac{2^{p+1}}{1-2^{p-2}}\|a_N\|_p^p\\
&\leq \frac{2^{p+1}}{1-2^{p-2}}\big(\frac{2^{p-1}}{2^{p-1}-1}\big)^p \|x_N\|_p^p.
\end{align*}
The desired inequality is achieved.
\end{proof}

%%%%%%%%%
%%%%%%%%%
 \section{Noncommutative Gundy's decomposition and differential subordination}
 
 In this section, we   present a new Gundy type decomposition that is well-suited for  dealing  with     weakly differentially subordinate martingales. We recall that  Gundy's decomposition  for noncommutative martingales was first considered in \cite[Theorem~2.1]{PR}. However, as we will see  below, the version given there does not easily fit with the notion of weak differential subordination.  We now state our new version:
\begin{theorem}\label{Gundy} 
Let  $x=(x_n)_{n\geq 1}$ be a self-adjoint  $L_1$-bounded martingale and $y$ is a self-adjoint martingale that is weakly differentially subordinate to $x$.  For any given  positive real number $\lambda$, there exist  four martingales $\a$, $\b$, $\g$, and $\upsilon$ satisfying 
the following properties:
\begin{enumerate}[{\rm (i)}]
\item $y=\a +\b +\g +\upsilon$;
\item  the martingale $\a$ satisfies:
$\|\a\|_2^2 \leq 2\lambda \| x\|_1$;
\item the martingale $\b$ satisfies:
\[
\sum_{n\geq 1} \|d\b_n\|_1 \leq 4 \|x\|_1;
\]
\item $\g$ and $\upsilon$ are $L_1$-martingales with:
\[
\max \Big\{ \lambda \T \Big( \bigvee_{n \ge 1} \mathrm{supp}
|d\g_n| \Big), \, \lambda \T \Big( \bigvee_{n \ge 1}
\mathrm{supp} \, |d\upsilon_n^*| \Big) \Big\} \leq
\|x\|_1.
\]
\end{enumerate}
\end{theorem}

We should point out here that according to \cite{Jiao-OS-Wu}, the martingale $y$ is not necessarily $L_1$-bounded and therefore the noncommutative Gundy's decomposition  from \cite{PR} does not apply directly to the martingale $y$.

\begin{proof} We consider the Cuculescu projections $(q_n^{(\lambda)})_{n\geq 1}$ relative to the martingale $x$. Below, we simply write $(q_n)$ for $(q_n^{(\lambda)})$.
We define
the martingales $\a$, $\b$, $\g$, and $\upsilon$ as
follows:

\begin{equation*}\label{decomposition}
\begin{cases}
d\a_n &:= q_{n-1} dy_n q_n - \E_{n-1} (q_{n-1} dy_n
q_n); \\
 d\b_n &:= q_{n-1}dy_n(q_{n-1}- q_{n})  -\E_{n-1}\big(q_{n-1}dy_n (q_{n-1}- q_{n})\big); \\
 d\g_n&:= dy_n ({\bf 1}-q_{n-1});\\ 
d\upsilon_n &:=({\bf 1} -q_{n-1})dy_n q_{n-1}.
\end{cases} \tag{$\mathbf{G}_{\lambda}$}
\end{equation*}
Clearly, $d\alpha$, $d\beta$, $d\gamma$, and $d\upsilon$ are
martingale difference sequences and $y=\alpha +\beta +\gamma
+\upsilon$. It is important to note here that   the four martingales in the decomposition depend on both martingales $x$ and $y$ as the Cuculescu projections were taken relative to $x$.

By Proposition~\ref{Cuculescu},  item ${\rm (iv)}$ is clearly satisfied.
It remains to verify items ${\rm (ii)}$ and ${\rm (iii)}$.
We begin with item ${\rm (ii)}$. 
 Using the  inequality $\T(|a-\E(a)|^2) \leq \T(|a|^2)$ for any conditional expectation $\E$ and an operator $a\in L_2(\M)$, we have for every $n\geq 1$,
\begin{align*}
\| d\a_n\|_2^2 &\leq \|q_{n-1}dy_n q_{n}\|_2^2 \\
&=\T\big(  q_{n}dy_nq_{n-1}dy_n q_{n}\big) \\
&=\T\big(  q_{n}[q_{n-1}dy_nq_{n-1}dy_n q_{n-1}]q_{n}\big).
\end{align*}
Since $y$ is weakly differentially subordinate to $x$, we obtain that
\begin{equation}\label{a-norm}
\| d\a_n\|_2^2 \leq  \T\big(  q_{n}dx_nq_{n-1}dx_n q_{n}\big).
\end{equation}
By Lemma~\ref{L2-norm}(i), we have for every $N\geq 1$,
\begin{equation}\label{a-finite}
\big\|\a_N\big\|_2^2=\sum_{n=1}^N \| d\a_n\|_2^2 \leq  \|q_Nx_N q_N\|_2^2 + 2\lambda\T\big(({\bf 1}-q_N)|x_N|\big).
\end{equation}
Using the fact that $\|q_Nx_N q_N\|_\infty \leq \lambda$, we further  get
\begin{align*}
\sum_{n=1}^N \| d\a_n\|_2^2 &\leq  \lambda \|q_Nx_N q_N\|_1 + 2\lambda\T\big(({\bf 1}-q_N)|x_N|\big)\\
&\leq \lambda \|q_Nx_N^+ q_N\|_1 + \lambda \|q_Nx_N^- q_N\|_1 + 2\lambda\T\big(({\bf 1}-q_N)|x_N|\big)\\
&=\lambda \T\big( q_N |x_N| q_N) +2\lambda\T\big(({\bf 1}-q_N)|x_N|\big)\\
&\leq 2\lambda \|x_N\|_1.
\end{align*}
Taking $N \to \infty$, we conclude that $
\|\a\|_2^2 \leq  2\lambda\|x\|_1$.

For  the martingale $\b$, we observe first  that since $y$ is weakly  differentially subordinated to $x$, we have for every $n\geq 1$,
\begin{align*}
| q_{n-1}dy_n(q_{n-1}- q_{n})|^2 &=
(q_{n-1}- q_{n}) dy_n q_{n-1}dy_n(q_{n-1}- q_{n})\\
&=(q_{n-1}- q_{n})[ q_{n-1}dy_n q_{n-1}dy_nq_{n-1}](q_{n-1}- q_{n})\\
&\leq (q_{n-1}- q_{n})[ q_{n-1}dx_n q_{n-1}dx_nq_{n-1}](q_{n-1}- q_{n})\\
&=(q_{n-1}- q_{n}) dx_n q_{n-1}dx_n(q_{n-1}- q_{n})\\
&=|q_{n-1}dx_n(q_{n-1}- q_{n})|^2.
\end{align*}
This shows  in particular that for every $n\geq 1$, 
\[
| q_{n-1}dy_n(q_{n-1}- q_{n})|\leq  |q_{n-1}dx_n(q_{n-1}- q_{n})|.
\]
  Moreover, by commutation (Proposition~\ref{Cuculescu}(ii)), the following equality holds:
  \[
  q_{n-1}dx_n(q_{n-1}- q_{n})=
(q_{n-1}- q_{n}) x_n (q_{n-1}- q_{n})
-q_{n-1}x_{n-1} (q_{n-1}- q_{n}). \]
We then  have the  following estimates for the $L_1$-norms:
\begin{align*}
\|q_{n-1}dy_n(q_{n-1}- q_{n})\|_1 &\leq 
\|q_{n-1}dx_n(q_{n-1}- q_{n})\|_1\\
&\leq \|(q_{n-1}- q_{n}) x_n (q_{n-1}- q_{n})\|_1 + \|q_{n-1}x_{n-1}q_{n-1} (q_{n-1}- q_{n})\|_1.
\end{align*}
Fix $N\geq 1$ and $n \leq N$. As $\|q_{n-1}x_{n-1}q_{n-1}\|_\infty \leq \lambda$, we get that,
\begin{align*}
\|q_{n-1}dy_n(q_{n-1}- q_{n})\|_1  &\leq \|(q_{n-1}- q_{n}) \E_n(x_N) (q_{n-1}- q_{n})\|_1  +\lambda \T\big(q_{n-1}- q_{n}\big) \\
&\leq \|(q_{n-1}- q_{n}) x_N (q_{n-1}- q_{n})\|_1  +\lambda \T\big(q_{n-1}- q_{n}\big).
\end{align*}
It follows  that for every $N\geq 1$, 
\[
\label{b-norm}
\sum_{n=1}^N \|q_{n-1}dy_n(q_{n-1}- q_{n})\|_1
\leq  \T\big(({\bf 1}- q_{N}) |x_N|\big) + \lambda\T\big( {\bf 1}- q_{N} \big).
\]
Using the fact that conditional expectations are contractions in $L_1(\M)$, we deduce that  for every $N\geq 1$,
\begin{equation}\label{b-norm}
\sum_{n= 1}^N \|d\b_n\|_1\leq 4\T\big(({\bf 1}- q_{N}) |x_N|\big).
\end{equation}
Taking the limit as $N\to \infty$,   we conclude that
$
\sum_{n\geq 1} \|d\b_n\|_1\leq 4\|x\|_1
$ as claimed.
\end{proof}

We observe that since  $y$ is an $L_1$-martingale that is not necessarily bounded,  the martingales $\g$ and $\upsilon$ are not necessarily $L_1$-bounded. 
However, the interest on $\g$ and $\upsilon$  is  only on their respective support projections. We also note that since $y$  is self-adjoint, we may take by symmetrization that $\a$ and $\b$ are self-adjoint martingales and $\g^*=\upsilon$.

\medskip

Even for the case $y=x$, the decomposition considered  here is different from the ones presented in \cite[Corollary~2.9]{PR}. It turns out that this new set up provides much better constants  on the norms of $\a$ and $\b$ than the decomposition previously considered in \cite{PR}. In fact, the constants obtained in items ${\rm (ii)}$ and ${\rm (iii)}$  are identical to those from the original decomposition for the classical case using one stoping time provided by Burkholder in \cite[Theorem~4.1]{Bu1}.  We also note that another  version of Gundy's decomposition was   considered in \cite{Jiao-OS-Wu} to accommodate the differential subordination but the version given there  was not made up of martingales and much less intuitive than Theorem~\ref{Gundy}.   

\medskip

We conclude this section  with  a natural  application of Theorem~\ref{Gundy}. We provide a very concise  alternative 
 proof of the weak type $(1,1)$ inequality for weak differential subordination  from \cite{Jiao-OS-Wu} in the spirit of the proof of weak type $(1,1)$ boundedness of  noncommutative martingale transforms from \cite[Theorem~3.1]{PR}. This approach also produces a better constant.
 
\begin{theorem}[{\cite[Theorem~4.1]{Jiao-OS-Wu}}]
Let $x$ be a self-adjoint $L_1$-bounded martingale and $y$ is a self-adjoint martingale that is  weakly differentially subordinate to $x$. Then for every $N\geq 1$,
\[
\big\|y_N\big\|_{1,\infty} \leq c \big\|x_N\big\|_1
\]
where $\displaystyle{c=2 +2B^2 + \frac{4B}{B-1} }$ for $B>1$.
\end{theorem}
 
 \begin{proof} By the definition of  the weak-$L_1$ norm, we need to verify that  for every $\lambda>0$, we have
 \[
 \lambda \T\big( \ch_{(\lambda,\infty)}(|y_N|) \big) \leq c\big\|x_N\big\|_1.
 \]
 By homogeneity,  it is enough to establish the inequality for $\lambda=1$. Consider the decomposition of $y=\a +\b + \g+\upsilon$ according to Theorem~\ref{Gundy} with $\lambda=1$. Fix $B>1$ and set $\eta=1/B$.

 From Lemma~\ref{Quasi-Triangle}, we have for every $0<\delta<1$,
 \begin{align*}
 \T\big(\ch_{(1,\infty)}(|y_N|)\big) &\leq  \T\big(\ch_{(\eta\delta,\infty)}(|\a_N|)\big) +
\T\big(\ch_{((1-\eta)\delta,\infty)}(|\b_N|)\big) \\
&+ \T\big(\ch_{(\eta(1-\delta),\infty)}(|\g_N|)\big) 
+\T\big(\ch_{((1-\eta)(1-\delta),\infty)}(|\upsilon_N|)\big) \\
&= I + II + III + IV.
\end{align*}
For the first term $I$, we use Chebychev's inequality to deduce:
\begin{equation} \label{I}
I \leq \frac{1}{\eta^2\delta^2} \|\a_N\|_2^2 \leq \frac{2}{\eta^2\delta^2} \big\| x_N\big\|_1.
\end{equation}
For the second term  $II$, we proceed similarly:
\begin{equation} \label{II}
\begin{split}
II &\leq  \frac{1}{(1-\eta)\delta} \big\|\b_N\big\|_1\\
&\leq  \frac{1}{(1-\eta)\delta}\sum_{n=1}^N \big\|d\b_n\big\|_1 \leq \frac{4}{(1-\eta)\delta}\big\|x_N\big\|_1.
\end{split}
\end{equation}
For $III$, we note that $|\g_N|=|\sum_{n=1}^N d\gamma_n|$ is
supported by the projection 
$\bigvee_{n=1}^N \mathrm{supp}
|d\gamma_n|$
and therefore  $\ch_{(\eta(1-\delta),\infty)}(|\g_N|)$ is a subprojection of $\bigvee_{n=1}^N \mathrm{supp}
|d\gamma_n|$.
 With this observation, it follows that
\begin{equation} \label{III}
III \leq  \T\Big(\bigvee_{n=1}^N \mathrm{supp}
|d\gamma_n|\Big) \leq  \big\|x_N\big\|_1.
\end{equation}
The last term $IV$ can be estimated similarly: 
\begin{equation} \label{IV}
\begin{split}
IV &=\T\big(\ch_{((1-\eta)(1-\delta),\infty)}(|\upsilon_N|)\big)\\
&=\T\big(\ch_{((1-\eta)(1-\delta),\infty)}(|\upsilon_N^*|)\big)\\
&\leq \T\Big(\bigvee_{n=1}^N \mathrm{supp}
|d\upsilon_n^*|\Big) \leq  \big\|x_N\big\|_1.
\end{split}
\end{equation}
Combining estimates  from
(\ref{I}, \ref{II}, \ref{III}, \ref{IV}) and taking $\delta \to 1$, we conclude that
\begin{align*}
\T\big(\ch_{(1,\infty)}(|y_N|)\big) &\leq \big( \frac{2}{\eta^2} +\frac{4}{1-\eta} +2\big) \big\|x_N\big\|_1\\
&=\big(2 +2B^2 + \frac{4B}{B-1}\big)\big\|x_N \big\|_1.
\end{align*}
The proof is complete.
 \end{proof}
\begin{remark}
Taking  for instance $B=1.75$ yields $c=17.458\overline{3}$ which is better than the constant $36$  in \cite[Theorem~4.1]{Jiao-OS-Wu}.
\end{remark}

%%%%%%%%%%%%%%%%%%%%%%%

\section{Differential subordinations and martingale Hardy spaces}

In this section, 
we present the primary objective of the paper. That is, to relate the  notion of weak  differential subordination with column and row square functions.  We organize  the section into three subsections. The first one is dedicated to  weak-type  results, the second deals with the corresponding strong type variants when $1<p<2$, while  in the last subsection, we present  noncommutative analogues of sharp constant result due to Wang \cite{Wang} which we use to
compare the various strong type results from the second subsection.

\subsection{Weak-type $(1,1)$ inequalities}\label{weak-type}

In this subsection, we provide two weak-type $(1,1)$ inequalities involving weak differential subordinations. Recall that as was discovered in \cite{Jiao-OS-Wu}, this version of noncommutative differential subordination is  needed when dealing with weak-type $(1,1)$ and strong type  for $1<p<2$. The inequalities considered below are motivated by \cite[Theorem~4.1]{Jiao-OS-Wu}. They are  in the spirit of the  weak-type  analogues  of the noncommutative Burkholder and the weak-type analogue of noncommutative Burkholder-Gundy from \cite{Ran18, Ran21}.

\begin{theorem}\label{main-weak}
Suppose that $x$ is a  self-adjoint $L_2$-martingale and $y$ is a  self-adjoint martingale that  is weakly  differentially subordinate to $x$. Then
 there exist
three adapted sequences $\eta=(\eta_n)_{n\geq 1}$, $\zeta=(\zeta_n)_{n\geq 1}$, and $\xi=(\xi_n)_{n\geq 1}$ such
that
$dy=\eta + \zeta +\xi$ and satisfy the weak-type estimate:
\[
\big\| \eta\big\|_{L_{1, \infty}(\M
\overline{\otimes} \ell_\infty)} + \big\| \zeta \big\|_{L_{1,\infty}^{\rm cond}(\M;\ell_2^c)} + 
\big\| \xi \big\|_{L_{1,\infty}^{\rm cond}(\M;\ell_2^r)}
\leq K \big\| x \big\|_1.
\]
\end{theorem}

\bigskip

Before we proceed with the proof, let us point out  that  as already mentioned in the introduction section, the corresponding result for classical martingales  can be  easily deduced  from previously known inequalities and therefore did not attract any   interest in general. Indeed, let $(\O,\Sigma,\mathbb{P})$ be a probability space  and $f=(f_n)_{n\geq 1}$ be a bounded martingale in $L_1(\O)$ (with respect to a  given filtration of  $\sigma$-subalgebras of $\Sigma$). Assume that $g=(g_n)_{n\geq 1}$  is a martingale that is differentially subordinate to $f$ in the sense of \cite{Burkholder-diff}. Following the classical Davis decomposition, define for $n\geq 1$,
\begin{equation*}
\begin{cases}
u_n &:=dg_n\ch_{\{f_n^*<2f_{n-1}^*\}} \\
v_n &:=dg_n\ch_{\{f_n^* \geq 2f_{n-1}^*\}}
\end{cases}
\end{equation*}
where $(f_n^*)_{n\geq 1}$ denotes the sequence of maximal functions. Then $dg_n=u_n +v_n$ and by  differential subordination, $|u_n| \leq |df_n|\ch_{\{f_n^*<2f_{n-1}^*\}}$ and $|v_n| \leq |df_n|\ch_{\{f_n^*\geq 2f_{n-1}^*\}}$. It then follows immediately from \cite[Corollary~C]{parcet} that:
\[
\Big\| \big( \sum_{n\geq 1}\mathbb{E}_{n-1} (|u_n|^2) \big)^{1/2}\Big\|_{1,\infty} + \sup_{\lambda>0}\Big\{ \lambda \sum_{n\geq 1} \mathbb{P}[ |v_n| \geq \lambda] \Big\}\leq C\big\|f\big\|_1.
\]
This clearly shows that Theorem~\ref{main-weak} is satisfied for classical martingales.

\bigskip

The noncommutative situation is very different. First, the observation above for the classical case suggests  that the natural approach  should follow the spirit of the decomposition considered in \cite{Ran21}. However, as we will see below, this does not lead to simple comparisons as in the classical case described above.

Below, we consider the family of projections $\{p_{i,n}\}_{i,n}$ constructed from the martingale $x$ as described in the previous section.
First, set  $\eta_1=y_1$, $\zeta_1=0$, and $ \xi_1=0$. For  $n\geq 2$, we define:
\begin{equation}\label{main-decom}
\begin{cases}
\eta_n &:=
\displaystyle{ \sum_{0\leq j<i}
(p_{i,n}-p_{i,n-1}p_{i,n}) dy_n p_{j,n-1}}; \\
\zeta_n &:=
\displaystyle{ \sum_{0\leq i \leq j} p_{i,n} dy_n
p_{j,n-1}}; \\
 \xi_n &:= 
    \displaystyle{\sum_{0\leq j<i}
p_{i,n-1}p_{i,n} dy_n p_{j,n-1}}.
\end{cases}
\end{equation}

Observe that by assumption, $y$ is  also an $L_2$-martingale. Since for every $n\geq 2$,  $dy_n \in L_2(\M)$, it is not difficult to verify  that the three sequences are well-defined and belong to $L_2(\M)$. This fact is  essential in computing their respective  norms in  conditioned  spaces. It is also worth pointing out here that  all three sequences depend  on both martingales $x$ and $y$ as the  family of projections $(p_{i,n})_{i,n}$ was derived from $x$.
Clearly, $\eta$, $\zeta$, and $\xi$ are adapted sequences and $dy=\eta +\zeta +\xi$.
We emphasize  that although $y$ is weakly differentially subordinate to $x$, that information does not carry over through the decomposition. In fact, none of the three adapted sequences consists of self-adjoint operators. In particular, for any $n\geq 2$, no immediate comparison can be made between say $|\zeta_n|$ (or $|\zeta_n^*|$) and $|dx_n|$.  Similar remark can be made with $\eta$ and $\xi$. 

\medskip

Our aim is to carefully analyze distribution functions of  appropriate operators relative to the three  sequences $\eta$, $\zeta$, and $\xi$. Without loss of generality, we may assume that both $x$ and $y$ are finite $L_2$-martingales. That is,   $x=(x_n)_{1\leq n\leq N}$ for some fixed $N$. Consequently, all three sequences in \eqref{main-decom} are finite sequences.

In this case,  $\sigma_c(\zeta)$  is a well-defined operator since the $\zeta_n$'s belong to $L_2(\M)$. By  the definition of column conditioned spaces, $\| \zeta \|_{L_{1,\infty}^{\rm cond}(\M;\ell_2^c)} =\| \sigma_c(\zeta)\|_{1,\infty}$. Similarly,
$\| \xi \|_{L_{1,\infty}^{\rm cond}(\M;\ell_2^r)} =\| \sigma_c(\xi^*)\|_{1,\infty}$.
On the other hand, we will view the sequence $\eta=(\eta_n)_{n=1}^N$  as an operator affiliated with  the von Neumann algebra $\M\overline{\otimes} \ell_\infty$. For convenience, we write $\eta=\sum_{n=1}^N \eta_n \otimes e_n$ where $(e_n)_{n\geq1}$ denotes the unit vector basis of $\ell_\infty$. Below,  $\Tr$  denotes  the natural trace on $\M\overline{\otimes} \ell_\infty$. That is, for $a=\sum_{n\geq 1} a_n \otimes e_n $, the trace of $a$ is given by 
$\Tr(a) =\sum_{n\geq 1} \T(a_n)$.

%%%%%%%
Thus, according  to the definition of the weak-$L_1$-norm,  proving Theorem~\ref{main-weak} amounts to getting the right estimates for the distribution functions of $\sigma_c(\zeta)$ and $\sigma_c(\xi^*)$ as $\T$-measurable operators and the distribution function of $\eta$ as $\Tr$-measurable operator. We record  these estimates in Proposition~\ref{distribution-c} and Proposition~\ref{distribution-d} below as they will also be needed in the next two subsections. We start with the column and row parts.

\begin{proposition}\label{distribution-c}
For every $k\geq 0$, the following two inequalities hold:
\begin{equation*}
\T\Big( \ch_{(2^k, \infty)}\big( \sigma_c(\zeta) \big) \Big) \leq 2^{-2k} \big\|e_{k,N }x_N e_{k,N}\big\|_2^2  +  8\ . \ 2^{-k}\T\big(({\bf 1} -e_{k,N})|x_N|\big)
\end{equation*}
and
\begin{equation*}
\T\Big( \ch_{(2^k, \infty)}\big( \sigma_c(\xi^*) \big) \Big) \leq 2^{-2k} \big\|e_{k,N }x_N e_{k,N}\big\|_2^2  +  8\ . \ 2^{-k}\T\big(({\bf 1} -e_{k,N})|x_N|\big).
\end{equation*}

\end{proposition}

%%%%%%

For the proof,
we  begin with a crucial lemma  that provides the decisive step  that transforms the norms of conditioned square functions of the sequences $\zeta$ and $\xi$ into square functions of some truncations of $x$. 

\begin{lemma}\label{transfer} For every $k\geq 0$ and $n\geq 2$, the following two inequalities hold:
\[
\T\Big(e_{k,{n-1}} \E_{n-1}(|\zeta_n|^2) e_{k,n-1} \Big) \leq \T\Big( e_{k,n}dx_n e_{k,n-1}dx_n e_{k,n}\Big)\] 
and 
\[
\T\Big(e_{k,{n-1}} \E_{n-1}(|\xi^*_n|^2) e_{k,n-1} \Big) \leq \T\Big( e_{k,n}dx_n e_{k,n-1}dx_n e_{k,n}\Big).\] 
\end{lemma}
\begin{proof} 
A straightforward computation  shows that for every $n\geq 2$,
\[
\E_{n-1}(|\zeta_n|^2) =\sum_{l,j\geq 0} \sum_{i\leq l\wedge j} p_{l,n-1}\E_{n-1}
[ dy_n p_{i,n} dy_n] p_{j,n-1}.\]
Fix $k\geq 0$. Since $e_{k,n-1}=\sum_{j=0}^k p_{j,n-1}$, we have:
\[
e_{k,{n-1}} \E_{n-1}(|\zeta_n|^2) e_{k,n-1}=\sum_{0\leq l,j\leq k} \sum_{i\leq l\wedge j} p_{l,n-1}\E_{n-1}
[ dy_n p_{i,n} dy_n] p_{j,n-1}.
\]
Taking the trace and using the fact that conditional expectations are trace invariant,
\begin{align*}
\T\Big(e_{k,{n-1}} \E_{n-1}(|\zeta_n|^2) e_{k,n-1} \Big) &=
\sum_{0\leq j\leq k} \sum_{0\leq i\leq j} \T\big(p_{j,n-1}\E_{n-1}
[ dy_n p_{i,n} dy_n] p_{j,n-1}\big)\\
&=\sum_{0\leq j\leq k}  \T\big(p_{j,n-1} dy_n e_{j,n} dy_n p_{j,n-1} \big).
\end{align*}
Since for $0\leq j\leq k$, $e_{j,n} \leq e_{k,n} \leq e_{k,n-1}$, we get
\begin{align*}
\T\Big(e_{k,{n-1}} \E_{n-1}(|\zeta_n|^2) e_{k,n-1} \Big)&\leq \sum_{0\leq j\leq k}  \T\big(p_{j,n-1} dy_n e_{k,n} dy_n p_{j,n-1} \big)\\
&=\sum_{0\leq j\leq k}  \T\big(p_{j,n-1} dy_n e_{k,n} dy_n \big)\\
&=  \T\big(e_{k,n-1} dy_n e_{k,n} e_{k,n-1}dy_n  \big)\\
&=  \T\big( e_{k,n} [e_{k,n-1}dy_n e_{k,n-1} dy_n e_{k,n-1}]\big)\\
&\leq   \T\big( e_{k,n} [e_{k,n-1}dx_n e_{k,n-1} dx_n e_{k,n-1}]\big)
\end{align*}
where in the last inequality we use the fact that $ e_{k,n-1} \in \M_{n-1}$ and $y$ is weakly differentially subordinate to $x$. This proves the first inequality. The second inequality is slightly more delicate. We begin from the identity:
\[
\E_{n-1}(|\xi_n^*|^2)
=\sum_{l,i\geq 0}\sum_{i<l\wedge j} p_{l,n-1}
\E_{n-1}[ p_{l,n}dy_n p_{j,n-1} dy_n p_{i,n}] p_{i,n-1}.
\]
 Performing the same computation  as in the first part leads to
\[
\T\Big(e_{k,{n-1}} \E_{n-1}(|\xi^*_n|^2) e_{k,n-1} \Big) 
=\sum_{0\leq i\leq k}  \T\big(p_{i,n-1} p_{i,n}dy_n e_{i,n-1} dy_n p_{i,n}p_{i,n-1} \big).
\]
Since $e_{i,n} \leq e_{k,n}$, we have
\[
\T\Big(e_{k,{n-1}} \E_{n-1}(|\xi^*_n|^2) e_{k,n-1} \Big) \leq \sum_{0\leq i\leq k}  \T\big(p_{i,n-1} p_{i,n}dy_n e_{k,n-1} dy_n p_{i,n}p_{i,n-1} \big).
\]
We also make the observation   that   $p_{i,n} \leq e_{i,n} \leq e_{k,n-1}$. This further leads to
\[
\T\Big(e_{k,{n-1}} \E_{n-1}(|\xi^*_n|^2) e_{k,n-1} \Big)
\leq  \sum_{0\leq i\leq k}  \T\big(p_{i,n-1} p_{i,n}[e_{k,n-1}dy_n e_{k,n-1} dy_n  e_{k,n-1}]p_{i,n}p_{i,n-1}\big). 
\]
Using the weak differential subordination assumption, we obtain
\begin{align*}
\T\Big(e_{k,{n-1}} \E_{n-1}(|\xi^*_n|^2) e_{k,n-1} \Big)
&\leq \sum_{0\leq i\leq k}  \T\big(p_{i,n-1} p_{i,n}[e_{k,n-1}dx_n e_{k,n-1} dx_n  e_{k,n-1}]p_{i,n}p_{i,n-1} \big)\\
&=\T\big([  \sum_{0\leq j\leq k} p_{i,n}p_{i,n-1} p_{i,n}]e_{k,n-1}dx_n e_{k,n-1} dx_n  e_{k,n-1}\big)\\
&\leq \T\big( \big[  \sum_{0\leq i\leq k}  p_{i,n}\big]e_{k,n-1}dx_n e_{k,n-1} dx_n  e_{k,n-1}\big)\\
&=\T\big(  e_{k,n}dx_n e_{k,n-1} dx_n  e_{k,n-1}\big).
\end{align*}
The proof of the lemma is complete.
\end{proof}

\begin{proof}[Proof of Proposition~\ref{distribution-c}] We present the proof for $\sigma_c(\zeta)$. The argument  for $\sigma_c(\xi^*)$ is identical.

Let  $k$  be a fixed positive integer. For simplicity, denote 
$\pi :=e_{k,N}$ and write \[
\sigma_c(\zeta)= \sigma_c(\zeta)\pi + \sigma_c(\zeta)({\bf 1}-\pi). 
\]
According to Lemma~\ref{Quasi-Triangle}, for any given $0<\alpha<1$,
\[
\T\Big( \ch_{(2^k, \infty)}\big( \sigma_c(\zeta) \big) \Big) \leq  \T\Big( \ch_{(\a 2^k, \infty)}\big( |\sigma_c(\zeta)\pi| \big) \Big) + \T\Big( \ch_{((1-\a)2^k, \infty)}\big( |\sigma_c(\zeta)({\bf 1}-\pi)| \big) \Big).
\]
Since $|\sigma_c(\zeta)({\bf 1}-\pi)|$ is supported by ${\bf 1}-\pi$, the projection 
$\ch_{((1-\a)2^k, \infty)}\big( |\sigma_c(\zeta)({\bf 1}-\pi)| \big)$ is a  subprojection of 
${\bf 1}-\pi$. Therefore,  we obtain that 
\[
\T\Big( \ch_{(2^k, \infty)}\big( \sigma_c(\zeta) \big) \Big) \leq  \T\Big( \ch_{(\a 2^k, \infty)}\big( |\sigma_c(\zeta)\pi| \big) \Big) +\T({\bf 1}-\pi).
\]
 Using Chebyshev's inequality and letting $\a \to 1$ on the first term of the right hand side, we get
\begin{equation}\label{col-2}
\begin{split}
\T\Big( \ch_{(2^k, \infty)}\big( \sigma_c(\zeta) \big) \Big) &\leq 2^{-2k}\T \big( \pi \sigma_c^2(\zeta)\pi \big) + \T({\bf 1}- e_{k,N})\\
&\leq 2^{-2k}\T \big( \pi \sigma_c^2(\zeta)\pi \big) + 2^{-k+1}\T\big(({\bf 1}- e_{k,N})|x_N|\big)
\end{split}
\end{equation}
where we use \eqref{e} on the second inequality.

Next, recall that  the sequence $(e_{k,m})$ is decreasing on $m$. In particular, $e_{k,m} \geq \pi$ for  every $1\leq m\leq N-1$.   With this fact, we  may express  the operator $\pi \sigma_c^2(\zeta) \pi$ as:
\begin{align*}
\pi \sigma_c^2(\zeta) \pi &=\pi\big( \sum_{n=1}^N \E_{n-1}(|\zeta_n|^2) \big)\pi\\
&=\pi\big( \sum_{n=2}^N e_{k,n-1}\E_{n-1}(|\zeta_n|^2)e_{k,n-1} \big)\pi.
\end{align*}
Taking  traces on both  sides leads to the inequality
\[
\T\big(\pi\sigma_c^2(\zeta) \pi) \leq \sum_{n=2}^N \T\big(e_{k,n-1}\E_{n-1}(|\zeta_n|^2)e_{k,n-1} \big).
\]
Now,  we appeal to  Lemma~\ref{transfer} to further get
\[
\T\big(\pi\sigma_c^2(\zeta) \pi) \leq \sum_{n=2}^N  \T\Big( e_{k,n}dx_n e_{k,n-1}dx_n e_{k,n}\Big).
\]
Therefore, we  see that \eqref{col-2}   further implies 
 \[
 \T\Big( \ch_{(2^k, \infty)}\big( \sigma_c(\zeta) \big) \Big) \leq  2^{-2k}\sum_{n=2}^N  \T\Big( e_{k,n}dx_n e_{k,n-1}dx_n e_{k,n}\Big) + 2^{-k+1}\T\big(({\bf 1}- e_{k,N})|x_N|\big).
 \]
We arrive at the desired estimate by incorporating the inequality  from Lemma~\ref{L2-norm}(iii)  on the first quantity of the right hand side.
\end{proof}

%%%%%

Now we deal with  the distribution function of $\eta$. The proof is  more involved.     Since the trace on $\M\overline{\otimes} \ell_\infty$  is not normalized, we need  to estimate the distribution function  for the full range of the positive real line. This is done in the next result.

\begin{proposition} \label{distribution-d} 
\begin{enumerate}[{\rm (i)}] 
\item For every $k\geq 1$, the following estimate holds:
\[
\Tr \big( \ch_{(2^k, \infty)}(|\eta|)\big) \leq 
 2^{-2k +2} \big\|e_{k,N} x_N e_{k,N}\big\|_2^2 +  28 \ . \  2^{-k} 
 \T\big(({\bf 1} -e_{k,N} )|x_N|\big) + \T\big( \ch_{(2^k,\infty)}(|x_1|)\big).
\]
\item For every   $0<\lambda\leq 1$, we have
\[
\Tr\big(\ch_{(\lambda,\infty)}(|\eta|)\big) \leq \T\big(\ch_{(\lambda, \infty)}(|x_1|)\big) +  \sum_{i\geq 0} \T({\bf 1}-e_{i,N}).
\]
\end{enumerate}
\end{proposition}

\begin{proof} $\bullet$ Assume that  $k\geq 1$.
Set $r_{i,n} :=r(p_{i,n}-p_{i,n-1}p_{i,n})$  and $l_{i,n}:=l(p_{i,n}-p_{i,n-1}p_{i,n})$ where $r(a)$ (resp. $l(a)$) denotes the right (resp. left) support projection of an operator $a$. It is clear that $r_{i,n} \leq p_{i,n}$. In particular, $(r_{i,n})_i$ is a pairwise disjoint sequence.
We claim that  for every $i\geq 1$,
\begin{equation}\label{left}
l_{i,n}  \leq e_{i-1,n-1}-e_{i-1,n}. 
\end{equation}
To verify this claim,  recall that  $p_{i,n}=e_{i,n}-e_{i-1,n}$ and $p_{i,n-1}= e_{i,n-1} -e_{i-1,n-1}$. Then 
\begin{align*}
p_{i,n}-p_{i,n-1}p_{i,n} &=e_{i,n}-e_{i-1,n}- ( e_{i,n-1} -e_{i-1,n-1})(e_{i,n}-e_{i-1,n})\\
&=e_{i,n}-e_{i-1,n}-e_{i,n-1}e_{i,n} +e_{i,n-1}e_{i-1,n}  +e_{i-1,n-1}e_{i,n} -e_{i-1,n-1}e_{i-1,n}. 
\end{align*}
Since $(e_{i,n})$ is decreasing on $n$ and increasing on $i$, we have $e_{i,n-1}e_{i,n}=e_{i,n}$, $e_{i-1,n-1}e_{i-1,n}=e_{i-1,n}$, and $e_{i,n-1}e_{i-1,n}=e_{i-1,n}$
where the last equality follows from $e_{i,n-1}\geq e_{i-1,n-1}\geq e_{i-1,n}$. It follows that the first four terms add up to zero and  therefore
\[
p_{i,n}-p_{i,n-1}p_{i,n} =e_{i-1,n-1}e_{i,n} -e_{i-1,n}. 
\]
Clearly, $e_{i-1,n-1}(p_{i,n}-p_{i,n-1}p_{i,n} )=p_{i,n}-p_{i,n-1}p_{i,n}$. Also, we have 
$e_{i-1,n}(p_{i,n}-p_{i,n-1}p_{i,n} )=e_{i-1,n}e_{i-1,n-1}e_{i,n} -e_{i-1,n}=0$. It follows that  $(e_{i-1,n-1}-e_{i-1,n})(p_{i,n}-p_{i,n-1}p_{i,n})=p_{i,n}-p_{i,n-1}p_{i,n}$ and therefore we have 
\eqref{left} from the definition of left support projection.

Now we proceed with  the proof.
To make the argument more transparent, we introduce the following notations:
\begin{align*}
U &:= {\bf 1} \otimes  e_1 + \sum_{n=2}^N \big(\sum_{0<i} p_{i,n}-p_{i,n-1}p_{i,n} \big) \otimes e_n;\\
V&:=y_1 \otimes e_1 + \sum_{n=2}^N \big(\sum_{0\leq j < i} r_{i,n} dy_n p_{j,n-1} \big) \otimes e_n;\\
\Pi_0 &:= {\bf 1} \otimes e_1 + \sum_{n=2}^N \big(\sum_{i > 0} r_{i,n} \big)  \otimes e_n;\\
\Pi_k &:= \chi_{(2^k, \infty)}(|x_1|) \otimes e_1 + \sum_{n=2}^N \big(\sum_{i\geq k+1} r_{i,n} \big) \otimes e_n.
\end{align*}
One can easily check  that $\eta=U \Pi_0 V$.  The operators  $\Pi_0$ and $\Pi_k$  are projections in $\M \overline{\otimes}\ell_\infty$  satisfying $ \Pi_k \leq \Pi_0$. We will need the following two properties:
\begin{enumerate}
\item $\Tr(\Pi_k) \leq \T(\ch_{(2^k,\infty)}(|x_1|)) + \sum_{i\geq k} \T \big({\bf 1}-e_{i,N}\big)$;
\item $\sup_{n,k} \big\| \sum_{i\geq k} p_{i,n} -p_{i,n-1} p_{i,n}\big\|_\infty \leq 2$.
\end{enumerate}
A verification of item $(2)$ can be found in \cite[Lemma~3.6]{Ran21}. On the other hand, we have from the definition that, 
\begin{align*}
\Tr(\Pi_k)&=\T(\ch_{(2^k,\infty)}(|x_1|)) + \sum_{n=2}^N \sum_{i\geq k+1} \T(r_{i,n})\\
&=\T(\ch_{(2^k,\infty)}(|x_1|)) + \sum_{n=2}^N \sum_{i\geq k+1} \T(l_{i,n}).
\end{align*}
Since $l_{i,n} \leq e_{i-1,n-1}-e_{i-1,n}$, a fortiori, we have
\begin{align*}
\Tr(\Pi_k)&\leq \T(\ch_{(2^k,\infty)}(|x_1|)) + \sum_{n=2}^N \sum_{i\geq k+1} \T(e_{i-1,n-1}-e_{i-1,n})\\
&\leq \T(\ch_{(2^k,\infty)}(|x_1|)) +  \sum_{i\geq k+1} \T({\bf 1}-e_{i-1,N}).
\end{align*}
Thus, we have deduced property~$(1)$.  That is,
\begin{equation}\label{big-proj}
\Tr(\Pi_k) \leq  \T(\ch_{(2^k,\infty)}(|x_1|)) +  \sum_{i\geq k} \T({\bf 1}-e_{i,N}).
\end{equation}
Next, we  write $\eta= U(\Pi_0 -\Pi_k)V + U\Pi_kV$ and fix $0<\a<1$. By Lemma~\ref{Quasi-Triangle}, we have
\[
\Tr\big(\ch_{(2^k,\infty)}(|\eta|)\big) \leq \Tr\big(\ch_{(\a2^k,\infty)}(|U(\Pi_0-\Pi_k)V|)\big) + \Tr\big(\ch_{((1-\a)2^k,\infty)}(|U\Pi_kV|)\big).
\]
We observe that by definition, the support projection of  $|U\Pi_kV|$ is a subprojection of $r(\Pi_kV)$. But since $r(\Pi_k V)$ is equivalent to $l(\Pi_kV)$ (see for instance \cite[p.~304]{TAK}) and $l(\Pi_kV)$ is clearly a subprojection of $\Pi_k$,
we see that $\ch_{((1-\a)2^k,\infty)}(|U\Pi_kV|)$ is equivalent to a subprojection of $ \Pi_k$ in $\M \overline{\otimes} \ell_\infty$.
Therefore, $\Tr(\ch_{((1-\a)2^k,\infty)}(|U\Pi_kV|)) \leq \Tr( \Pi_k)$.
 Applying Chebyshev inequality and taking $\a\to 1$ in the first term, we obtain that
\[
\Tr\big(\ch_{(2^k,\infty)}(|\eta|)\big) \leq  2^{-2k}\big\|U(\Pi_0-\Pi_k)V\big\|_2^2  + \Tr\big(\Pi_k\big).
\]
We proceed to estimate  the $L_2$-norm of $U(\Pi_0-\Pi_k)V$.
\begin{align*}
\big\|U(\Pi_0-\Pi_k)V\big\|_2^2  &=\big\|y_1\ch_{[0, 2^k]}(|x_1|)\big\|_2^2 +\sum_{n=2}^N
\big\| \sum_{0\leq j < i\leq k} (p_{i,n}-p_{i,n-1}p_{i,n} )dy_n p_{j,n-1}\big\|_2^2 \\
&=\big\|y_1\ch_{[0, 2^k]}(|x_1|)\big\|_2^2 +\sum_{n=2}^N
\big\| \big(\sum_{l\geq 0} p_{l,n}-p_{l,n-1}p_{l,n}\big)\sum_{0\leq j < i\leq k} p_{i,n}dy_n p_{j,n-1}\big\|_2^2\\
&=\big\|y_1\ch_{[0, 2^k]}(|x_1|)\big\|_2^2 +\sum_{n=2}^N
 \big\|\sum_{l\geq 0} p_{l,n}-p_{l,n-1}p_{l,n}\big\|_\infty^2 \big\|\sum_{0\leq j < i\leq k} p_{i,n}dy_n p_{j,n-1}\big\|_2^2.
\end{align*}
By the property~(2) stated above, $\|\sum_{l\geq 0} p_{l,n}-p_{l,n-1}p_{l,n}\|_\infty\leq 2$. Therefore,
\begin{align*}
\big\|U(\Pi_0-\Pi_k)V\big\|_2^2&\leq \big\|y_1\ch_{[0, 2^k]}(|x_1|)\big\|_2^2 +4\sum_{n=2}^N \big\|\sum_{0\leq j < i\leq k} p_{i,n}dy_n p_{j,n-1}\big\|_2^2\\
&=\big\|y_1\ch_{[0, 2^k]}(|x_1|)\big\|_2^2 +4\sum_{n=2}^N \big\|\sum_{1\leq  i\leq k} p_{i,n}dy_n e_{i-1,n-1}\big\|_2^2\\
&=\big\|y_1\ch_{[0, 2^k]}(|x_1|)\big\|_2^2 +4\sum_{n=2}^N \sum_{1\leq  i\leq k} \big\|p_{i,n}dy_n e_{i-1,n-1}\big\|_2^2.
\end{align*}
Since $e_{i-1,n-1} \leq e_{k,n-1}$ and  $(p_{i,n})_i$ are pairwise disjoint, we further get
\begin{align*}
\big\|U(\Pi_0-\Pi_k)V\big\|_2^2 &\leq \big\|y_1\ch_{[0, 2^k]}(|x_1|)\big\|_2^2 +4\sum_{n=2}^N \big\|e_{k,n}dy_n e_{k,n-1}\big\|_2^2\\
 &=\big\|y_1\ch_{[0, 2^k]}(|x_1|)\big\|_2^2 +4\sum_{n=2}^N \T\big(e_{k,n}[e_{k,n-1}dy_n e_{k,n-1}dy_n e_{k,n-1}] e_{k,n}\big).
\end{align*}
By assumption, $e_{k,n-1}dy_n e_{k,n-1}dy_n e_{k,n-1}\leq e_{k,n-1}dx_n e_{k,n-1}dx_n e_{k,n-1}$. Also $|y_1|^2 \leq |x_1|^2$. We  deduce  that 
\[
\big\|U(\Pi_0-\Pi_k)V\big\|_2^2\leq \big\|x_1\ch_{[0, 2^k]}(|x_1|)\big\|_2^2 +4\sum_{n=2}^N \T\big(e_{k,n}[e_{k,n-1}dx_n e_{k,n-1}dx_n e_{k,n-1}] e_{k,n}\big)
\]
which is equivalent to
\begin{equation}\label{big-L2}
\big\|U(\Pi_0-\Pi_k)V\big\|_2^2\leq \big\|x_1\ch_{[0, 2^k]}(|x_1|)\big\|_2^2 +4\sum_{n=2}^N \big\|e_{k,n}dx_n e_{k,n-1}\big\|_2^2.
\end{equation}
At this point,  we have from combining \eqref{big-proj} and \eqref{big-L2} that:
\begin{equation*}
\begin{split}
\Tr\big(\ch_{(2^k,\infty)}(|\eta|)\big) &\leq 2^{-2k +2} \sum_{n=2}^N \big\|e_{k,n}dx_n e_{k,n-1}\big\|_2^2+  \sum_{i\geq k} \T({\bf 1}-e_{i,N})\\
 &\ +
2^{-2k}\big\|x_1\ch_{[0, 2^k]}(|x_1|)\big\|_2^2 +
 \T(\ch_{(2^k,\infty)}(|x_1|)).
 \end{split}
 \end{equation*}
To simplify the right hand side,  we will incorporate the  third  term into the first term. To do this, we simply observe that at the first level of the Cuculescu construction, only one operator is being used and therefore all projections involved commute.  In particular,  we have $\ch_{[0,2^k]}(|x_1|)=q_1^{(2^k)}=e_{k,1}$. Therefore,
$\|x_1 \ch_{[0,2^k]}(|x_1|)\|_2^2 =\|e_{k,1}dx_1 e_{k,0} \|_2^2$.
 This allows us to write:
 \[
 \Tr\big(\ch_{(2^k,\infty)}(|\eta|)\big) \leq 2^{-2k +2} \sum_{n=1}^N \big\|e_{k,n}dx_n e_{k,n-1}\big\|_2^2+  \sum_{i\geq k} \T({\bf 1}-e_{i,N}) +
 \T(\ch_{(2^k,\infty)}(|x_1|)).
 \]
 By Lemma~\ref{L2-norm}(iii) and \eqref{e},
\begin{equation*}
\begin{split}
\Tr\big(\ch_{(2^k,\infty)}(|\eta|)\big)  &\leq 2^{-2k +2} \big\|e_{k,N} x_N e_{k,N}\big\|_2^2  +  24\ . \ 2^{-k}\T\big(({\bf 1}-e_{k,N})|x_N|\big) \\
 &+ \sum_{i\geq k} 2^{-i+1} \T\big(({\bf 1}-e_{i,N})|x_N|\big)  + \T(\ch_{(2^k,\infty)}(|x_1|)).
 \end{split}
 \end{equation*}
Since ${\bf 1}-e_{i,N}\leq {\bf 1}-e_{k,N}$ for $i\geq k$, we conclude that
\[
\Tr\big(\ch_{(2^k,\infty)}(|\eta|)\big)\leq 2^{-2k +2} \big\|e_{k,N} x_N e_{k,N}\big\|_2^2  +  28\ . \  2^{-k}\T\big(({\bf 1}-e_{k,N})|x_N|\big) +\T(\ch_{(2^k,\infty)}(|x_1|)).
\]
This completes the proof of the first item.
%%%%%

%%%
$\bullet$ Fix $0<\lambda\leq 1$.  Let $\Theta:= \eta-(y_1 \otimes e_1)$. It is clear from disjointness that $\ch_{(\lambda, \infty)}(|\eta|)=\ch_{(\lambda, \infty)}(|y_1|) \otimes e_1 +
\ch_{(\lambda,\infty)}(|\Theta|)$. Since  we clearly have by assumption that $|y_1| \leq |x_1|$, it follows that $\T(\ch_{(\lambda, \infty)}(|y_1|))\leq \T(\ch_{(\lambda, \infty)}(|x_1|))$. We claim that the inequality 
$\Tr(\ch_{(\lambda,\infty)}(|\Theta|)) \leq \sum_{i\geq 0}\T({\bf 1}-e_{i,N})$ holds. To verify this claim, we consider a slight modification of the projection $\Pi_0$ defined earlier. Set $\widetilde{\Pi}_0 :=\sum_{n=2}^N (\sum_{i > 0} r_{i,n}) \otimes e_n$.  We observe that
\[
\Theta= U\widetilde{\Pi}_0 V.
\]
Then, as before, the  support projection of $|\Theta|$ is a subprojection of $r(\widetilde{\Pi}_0 V)$
which in turn is equivalent to $l (\widetilde{\Pi}_0V)$. But $ l(\widetilde{\Pi}_0V)\leq \widetilde{\Pi}_0$. This shows that $\ch_{(\lambda,\infty)}(|\Theta|)$ is equivalent to a subprojection of $\widetilde{\Pi}_0$. This leads to
\[
\Tr\big(\ch_{(\lambda,\infty)}(|\Theta|)\big) \leq  \Tr\big(\widetilde{\Pi}_0\big)
=\sum_{n=2}^N \sum_{i>0} \T(r_{i,n}).
\]
We can deduce as before that 
$\sum_{n=2}^N \sum_{i>0} \T(l_{i,n})\leq 
\sum_{i>0} \T({\bf 1}- e_{i-1,N})= \sum_{i\geq 0} \T({\bf 1}-e_{i,N})$.
The proof  of the proposition is complete.
\end{proof}

%%%%%%%%%%%%%
%%%%%%%%%%%%%
We are now ready to provide the proof of the main theorem.

\begin{proof}[Proof of Theorem~\ref{main-weak}]
We will divide the proof into the column/row part and the diagonal part.

For the column part, we will show that there is an absolute constant $K$ so that for every $\lambda>0$,
\begin{equation}\label{column-inequality}
\lambda \T\Big( \ch_{(\lambda, \infty)}\big( \sigma_c(\zeta) \big) \Big) \leq K \|x\|_1.
\end{equation}
Since $\M$ is a finite von Neumann algebra, it suffices to establish \eqref{column-inequality} for $\lambda=2^k$ when $k\geq 0$ is arbitrary.
Note that $\|e_{k,N}x_N e_{k,N}\|_2^2 \leq 2^k \|x_N\|_1$. This  simple inequality
  along with Proposition~\ref{distribution-c}  clearly lead  to the conclusion that 
\begin{equation}\label{zeta-weak} 
2^k\T\Big( \ch_{(2^k, \infty)}\big( \sigma_c(\zeta) \big) \Big) \leq  9\|x\|_1.
\end{equation}
 The proof for the column part is complete. In light of the second inequality in Proposition~\ref{distribution-c}, the proof for the row part is identical. That is,
\begin{equation}\label{xi-weak} 
2^k\T\Big( \ch_{(2^k, \infty)}\big( \sigma_c(\xi^*) \big) \Big) \leq  9\|x\|_1.
\end{equation}

We  now provide the proof of the diagonal part. Proposition~\ref{distribution-d}(i) gives for  every $k\geq 1$ that
\[
\Tr\big(\ch_{(2^k,\infty)}(|\eta|)\big)  \leq  2^{-2k +2}( 2^{k} \|x_N\|_1) + 29\ . \ 2^{-k}\|x_N\|_1 = 33 \big(2^{-k}\|x_N\|_1\big).
\]
On the other hand,  for every $0<\lambda \leq 1$, Proposition~\ref{distribution-d}(ii) and \eqref{e} imply that
\[
\Tr\big(\ch_{(\lambda,\infty)}(|\eta|)\big) \leq \lambda^{-1}\|x_1\|_1 + 2\big(\sum_{i\geq 0}2^{-i} \|x_N\|_1\big) \leq 5\big(\lambda^{-1} \|x_N\|_1\big).
\]
Combining the previous two estimates, we can conclude that for every $l \in \mathbb{Z}$,
\begin{equation}\label{eta-weak}
2^l\Tr\big(\ch_{(2^l,\infty)}(|\eta|)\big) )\leq 33\|x_N\|_1.
\end{equation}
Thus,  getting  the desired estimate for the diagonal part. Finally, combining \eqref{zeta-weak}, \eqref{xi-weak}, and \eqref{eta-weak}, we conclude that 
\[
\big\| \zeta\big\|_{L_{1,\infty}^{\rm{cond}}(\M;\ell_2^c)} +\big\| \xi\big\|_{L_{1,\infty}^{\rm{cond}}(\M;\ell_2^r)} + \big\| \eta\big\|_{L_{1,\infty}^{\rm{cond}}(\M \overline{\otimes}\ell_\infty)} \leq 51 \big\| x\big\|_1.
\]
The proof is complete.
\end{proof}

\bigskip

Our next result is the companion of  Theorem~\ref{main-weak} for the case of square functions. 

\begin{theorem}\label{main-weak-S}
Suppose that $x$ is a  self-adjoint $L_2$-martingale and $y$ is a  self-adjoint martingale that  is weakly  differentially subordinate to $x$. Then
 there exist
two martingales $y^r$ and $y^c$ such
that
$y=y^c +y^r$ and satisfy the weak-type estimate:
\[
\big\| S_c(y^c)\big\|_{1, \infty} + \big\| S_r(y^r)\big\|_{1, \infty} 
\leq K \big\| x \big\|_1.
\]
\end{theorem}

We now proceed with the proof of Theorem~\ref{main-weak-S}. The martingales $y^c$ and $y^r$
are defined from their respective martingale difference sequences
as follows:
\begin{equation}\label{equation-square}
\begin{cases} dy_1^c &:=\displaystyle{ 
\sum_{0\leq i \leq j} p_{i,1} dy_1 p_{j,1}}; 
 \\ dy_n^c &:=\displaystyle{
\sum_{0\leq i \leq j} p_{i,n-1} dy_n p_{j,n-1}}  \quad
\text{for $n\geq 2$}; \\  dy_1^r &:=
\displaystyle{\sum_{0\leq j<i} p_{i,1} dy_1 p_{j,1}};
\\ dy_n^r &:= \displaystyle{\sum_{0\leq j<i}
p_{i,n-1} dy_n p_{j,n-1}} \quad \text{for $n\geq 2$}.
\end{cases}
\end{equation}
Clearly, $y=y^c +y^r$.  Since $x$ is an $L_2$-martingale, by the weak differential subordination assumption, $y$ is also an $L_2$-martingale.  By  the $L_2$-boundedness  of triangular truncations, it is also clear that $y^c$ and $y^r$ are $L_2$-martingales. 
However, since $y^c$ and $y^r$ are not self-adjoint martingales, neither of them can be weakly differentially subordinate to $x$.

\medskip

We will  present the argument for $\|S_c(y^c)\|_{1,\infty}$.  We emulate  the idea  of the proof of \cite[Theorem~3.2]{PR} taking advantage of  the  new version of Gundy's decomposition for  weakly differentially subordinate  martingales presented in the previous section.  Before we proceed with the proof, we need to  an intermediary result on triangular truncations.

Let $\mathsf{P} = \{ p_i \}_{i=1}^m$ be a finite sequence of
mutually disjoint projections in $\M$. We denote by $\mathcal{T}^{(\mathsf{P})}$ the  
\emph{triangular truncation} with respect to $\mathsf{P}$. That is, for any operator  $a\in L_0(\M,\T)$, we set:
\[ \mathcal{T}^{(\mathsf{P})} a =  \sum_{1\leq i\leq j\leq m} p_i a
p_j . \] 
We will make use  of the  following  property of countable set of   triangular truncations:
\begin{lemma}\label{Truncation1} If
$(\mathsf{P}_n)_{n\geq 1}$ is a family of finite sequences of
mutually disjoint  projections and $(a_n)_{n\geq 1}$  is a
sequence in $L_1(\M,\T)$, then
\[
 \Big\| \Big( \sum_{n\geq 1} \big| \mathcal{T}^{(\mathsf{P}_n)} a_n
\big|^2 \Big)^{1/2} \Big\|_{1,\infty} \leq 20\sqrt{2} \sum_{n\geq
1}\big \|a_n \big\|_1.
\]
\end{lemma}
A version of Lemma~\ref{Truncation1} for sequences of positive elements of $L_1(\M)$ appeared in \cite[Proposition~1.6]{Ran18} with constant equals to  $5\sqrt{2}$.
By the quasi-triangle inequality $\|w_1 +w_2\|_{1,\infty} \leq 2\|w_1\|_{1,\infty} + 2\|w_2\|_{1,\infty}$, one can easily verify that the lemma holds with constant equals to $10\sqrt{2}$ for  sequences of self-adjoint operators. Finally, Lemma~\ref{Truncation1} follows from  splitting the sequence $(a_n)_{n\geq 1} $ into real and imaginary parts.

%%%%
We are now ready to present the proof.
As in the conditioned case,  it suffices to  establish that there exists a constant $K$ such that  for  any arbitrary $k\geq 0$:

\begin{equation}\label{column-square}
2^k \T\Big( \ch_{(2^k, \infty)}\big( S_c(y^c) \big) \Big) \leq K \|x\|_1.
\end{equation}

Fix $k\geq 0$ and $N\geq 2$. For simplicity,  let $\pi=e_{k,N}$.  First, we write
\[
S_{c,N}(y^c)= S_{c,N}(y^c)\pi + S_{c,N}(y^c)({\bf 1}-\pi). 
\]
According to Lemma~\ref{Quasi-Triangle}, for any given $0< \delta<1$,
\begin{equation}\label{distribution-S1}
\begin{split}
\T\Big( \ch_{(2^k, \infty)}\big( S_{c,N}(y^c) \big) \Big) &\leq  \T\Big( \ch_{(\delta 2^k, \infty)}\big( |S_{c,N}(y^c)\pi| \big) \Big) + \T\Big( \ch_{((1-\delta)2^k, \infty)}\big( |S_{c,N}(y^c)({\bf 1}-\pi)| \big) \Big)\\
&\leq  \T\Big( \ch_{(\delta 2^k, \infty)}\big( |S_{c,N}(y^c)\pi| \big) \Big) + \T({\bf 1}-\pi).
\end{split}
\end{equation}
We note that $\pi S^2_{c,N}(y^c) \pi = \pi (S_{c,N}^{(k)}(y^c))^2 \pi$  where $S_{c,N}^{(k)}(y^c)$
denotes the following  truncated square
function 
\[
S_{c,N}^{(k)}(y^c) = \left( \Big| 
\sum_{0\leq i\leq j \leq k} p_{i,1} dy_1 p_{j,1} \Big|^2 + \sum_{n= 2}^N
\Big|  \sum_{0\leq i\leq j\leq k} p_{i,n-1} dy_n p_{j,n-1}
\Big|^2 \right)^{1/2}.
\]
We refer to the proof of \cite[Proposition~A]{Ran18} for this fact. Next, we consider the  decomposition  of $y$ according to Theorem~\ref{Gundy} using  the parameter $\lambda = 2^k$. This gives $y =
\alpha + \beta + \gamma + \upsilon$ as  described in (\ref{decomposition}).

As in \cite{PR}, for $n \geq 1$, set
\[
\mathsf{P}^{(k)}_n :=(p_{i,n})_{i=0}^k.
\]
Then with this notation, we see that
\[
S_{c,N}^{(k)}(y^c) = \Big( \big|
\cal{T}^{\mathsf{P}^{(k)}_1}(dy_1) \big|^2 + \sum_{n=2}^N
\big| \cal{T}^{\mathsf{P}^{(k)}_{n-1}}(dy_n) \big|^2 \Big)^{1/2}.
\]
We make the following crucial observation. Since $d\gamma_n$ is right supported by ${\bf 1}-q_{n-1}^{(2^k)}$. A fortiori, it is right supported by ${\bf 1}-e_{k,n-1}=\sum_{l\geq k+1} p_{k,n-1}$. This reveals that $\cal{T}^{\mathsf{P}^{(k)}_{n-1}}(d\gamma_n)=0$. Similarly, by using the left support projections, we see that $\cal{T}^{\mathsf{P}^{(k)}_{n-1}}(d\upsilon_n)=0$. Therefore,
\[
S_{c,N}^{(k)}(y^c) = \Big( \big|
\cal{T}^{\mathsf{P}^{(k)}_1}(d\a_1 +d\beta_1) \big|^2 + \sum_{n= 2}^N
\big| \cal{T}^{\mathsf{P}^{(k)}_{n-1}}(d\a_n + d\beta_n) \big|^2 \Big)^{1/2}.
\]
Using the elementary identity $|a+b|^2 \leq 2 |a|^2 +2|b|^2$, we have 
\[
(S_{c,N}^{(k)}(y^c))^2 \leq  2 (S_{c,N}^{(k)}(\a))^2 +2(S_{c,N}^{(k)}(\beta))^2
\]
where we use the notation

\[
S_{c,N}^{(k)}(\a) = \Big( \big|
\cal{T}^{\mathsf{P}^{(k)}_1}(d\a_1) \big|^2 + \sum_{n=2}^N
\big| \cal{T}^{\mathsf{P}^{(k)}_{n-1}}(d\a_n) \big|^2 \Big)^{1/2}
\]
and
\[
S_{c,N}^{(k)}(\beta) = \Big( \big|
\cal{T}^{\mathsf{P}^{(k)}_1}(d\beta_1) \big|^2 + \sum_{n=2}^N
\big| \cal{T}^{\mathsf{P}^{(k)}_{n-1}}( d\beta_n) \big|^2 \Big)^{1/2}.
\]
The preceding discussion and  Lemma~\ref{Quasi-Triangle} lead to:
\begin{align*}
\T\Big( \ch_{(\delta 2^k, \infty)}\big( |S_{c,N}(y^c)\pi| \big) \Big)&=\T\Big( \ch_{(\delta 2^k, \infty)}\big( |S_{c,N}^{(k)}(y^c)\pi| \big) \Big)\\
&\leq \T\Big( \ch_{(\delta^2 2^{2k}, \infty)}\big( |S^{(k)}_{c,N}(y^c)|^2 \big) \Big)\\
&\leq \T\Big( \ch_{(\delta^2 2^{2k-1}, \infty)}\big(2 |S^{(k)}_{c,N}(\a)|^2 \big) \Big) +\T\Big( \ch_{(\delta^2 2^{2k-1}, \infty)}\big(2 |S^{(k)}_{c,N}(\beta)|^2 \big) \Big)\\
&=\T\Big( \ch_{(\delta^2 2^{2k-2}, \infty)}\big( |S^{(k)}_{c,N}(\a)|^2 \big) \Big) +\T\Big( \ch_{(\delta^2 2^{2k-2}, \infty)}\big( |S^{(k)}_{c,N}(\beta)|^2 \big) \Big).
\end{align*}
Using Chebychev's inequality,  we further get:
\begin{align*}
\T\Big( \ch_{(\delta 2^k, \infty)}\big( |S_{c,N}(y^c)\pi| \big) \Big)&\leq  \delta^{-2} 2^{-2k+2} \| S^{(k)}_{c,N}(\a)\|_2^2 +\T\Big( \ch_{(\delta 2^{k-1}, \infty)}\big( S^{(k)}_{c,N}(\beta) \big) \Big)\\
&\leq \delta^{-2} 2^{-2k +2} \| S^{(k)}_{c,N}(\a)\|_2^2 +\delta^{-1}2^{-k+1} \|S^{(k)}_{c,N}(\beta)\|_{1,\infty}.
\end{align*}
Combining this last estimate with \eqref{distribution-S1} and taking $\delta \to 1$, we have
\[
\T\Big( \ch_{(2^k, \infty)}\big( S_{c,N}(y^c) \big) \Big) \leq 2^{-2k+2} \| S^{(k)}_{c,N}(\a)\|_2^2 +2^{-k+1} \|S^{(k)}_{c,N}(\beta)\|_{1,\infty} +\T({\bf 1}-e_{k,N}).
\]
Using the fact that triangular truncations are contractive projections in $L_2(\M)$ on the first term on the right hand side and Lemma~\ref{Truncation1} on the second term, we further get
\begin{equation}\label{distribution-Gundy}
\T\Big( \ch_{(2^k, \infty)}\big( S_{c,N}(y^c) \big) \Big) \leq 2^{-2k+2}\| \a_N\|_2^2 + (20\sqrt{2}) 2^{-k+1} \sum_{n=1}^N \|d\beta_n\|_1 + \T({\bf 1}-e_{k,N}).
\end{equation}
We can conclude from \eqref{a-finite}, \eqref{b-norm},  and \eqref{e} that 
\begin{align*}
\T\Big( \ch_{(2^k, \infty)}\big( S_{c,N}(y^c) \big) \Big)  &\leq 2^{-k +3}\| x_N\|_1 + (80\sqrt{2}) 2^{-k+1 } \|x_N\|_1+  2^{-k +1}\|x_N\|_1 \\
&=[10 +160\sqrt{2}] 2^{-k} \|x_N\|_1.
\end{align*}
Taking the limit with $N\to \infty$, we obtain  \eqref{column-square} with $K=10 +160\sqrt{2}$.
The proof is complete.
\qed

\begin{remark}\label{distribution-similar}
Inspecting the proof of Theorem~\ref{Gundy}, more specifically, the estimates on the $L_2$-norms of $(d\a_n)_{n\geq 1}$ and the $L_1$-norms of $(d\beta_n)_{n\geq 1}$, one can show  that the distribution function of $S_c(y^c)$  can be majorized with an upper bound similar to that of the distribution of  $\sigma_c(\zeta)$ from Proposition~\ref{distribution-c}. More precisely, we have for every $N\geq 1$,
\begin{equation}\label{dist-S}
\T\Big( \ch_{(2^k, \infty)}\big( S_{c,N}(y^c) \big) \Big) \leq 2^{-2k +2} \big\| e_{k,N}x_N e_{k,N}\big\|_2^2 + [160\sqrt{2} + 26] 2^{-k}\T\big( ({\bf 1}- e_{k,N} )|x_N|\big).
\end{equation}
The same  estimate  applies to  the distribution function of $S_{r,N}(y^r)$.
These are  nearly identical to  the estimates from Proposition~\ref{distribution-c} and will be useful in the extensions to convex functions and   the strong type $(p,p)$  case below. 
\end{remark}
We will sketch the argument for inequality \eqref{dist-S}. First, it follows from the proof of Lemma~\ref{L2-norm} and \eqref{a-norm} that  for $N\geq 1$,
\begin{equation}\label{a2-norm}
\|\a_N\|_2^2  \leq  \big\| e_{k,N}x_N e_{k,N}\big\|_2^2 + 6\ .\ 2^k \T\big( ({\bf 1}- e_{k,N} )|x_N|\big).
\end{equation}
Also since ${\bf 1}-q_N^{(2^k)} \leq {\bf 1}-e_{k,N}$, it follows from \eqref{b-norm} that 
\begin{equation}\label{b2-norm}
\sum_{n=1}^N \|d\b_n\|_1 \leq 4 \T\big( ({\bf 1}-e_{k,N})|x_N|\big).
\end{equation}
Incorporating \eqref{a2-norm} and \eqref{b2-norm}  into the inequality \eqref{distribution-Gundy}, we obtain 
\begin{equation*}
\begin{split}
\T\Big( \ch_{(2^k, \infty)}\big( S_{c}(y^c) \big) \Big) &\leq 2^{-2k +2} \Big[\big\| e_{k,N}x_N e_{k,N}\big\|_2^2 + 6\ .\ 2^k \T\big( ({\bf 1}- e_{k,N} )|x_N|\big)\Big] \\
&+ (20\sqrt{2})2^{-k+1} \Big[4 \T\big( ({\bf 1}- e_{k,N} )|x_N|\big)\Big]   + \T\big( {\bf 1}- e_{k,N} \big)\\
&\leq  2^{-2k +2} \big\| e_{k,N}x_N e_{k,N}\big\|_2^2 +\big[ 160\sqrt{2} +26\big]2^{-k} \T\big( ({\bf 1}- e_{k,N} )|x_N|\big)
\end{split}
\end{equation*}
where in the second inequality we use \eqref{e}. \qed

%%%%%%%%
%%%%%%%%

\subsection{Strong type $(p,p)$  inequalities for $1<p<2$}
In this subsection, we consider  the corresponding  strong-type $(p,p)$ for weakly differentially subordinate martingales. Our first result deals with conditioned Hardy spaces which is  the strong type version of Theorem~\ref{main-weak}.

\begin{theorem}\label{strong-type-1} 
Suppose that $x$ is a  self-adjoint $L_2$-martingale and $y$ is a  self-adjoint martingale that  is weakly  differentially subordinate to $x$. Then
 there exist
three martingales $y^d$, $y^c$, and $y^r$ (depending only on  $x$ and $y$) such
that
$y=y^d + y^c +y^r$ and  for every $1<p<2$,
\[
\big\| y^d\big\|_{\h_p^d} + \big\| y^c\big\|_{\h_p^c}  +\big\| y^r\big\|_{\h_p^r} \leq c_p \big\| x \big\|_p
\]
where 
\[
c_p =\frac{2^{p+1}}{(2^{p-1}-1)} \Big[  \Big(8 +\frac{6}{1-2^{p-2}}\Big)^{1/p} + \Big(77 +\frac{24}{1-2^{p-2}}\Big)^{1/p}\Big].
\]
In particular, $c_p=O((p-1)^{-1})$ when $p \to 1$.
\end{theorem}

  We divide  the proof into two steps:

{\bf Step~1.} 
We begin by describing a  concrete decomposition of $y$. Consider  the martingales $w^d$, $w^c$,  and $w^r$ whose martingale difference sequences  are for $n \geq 1$, $dw_n^d :=\eta_n -\E_{n-1}(\eta_n)$, $dw_n^c :=\zeta_n -\E_{n-1}(\zeta_n)$, and $dw_n^r :=\xi_n -\E_{n-1}(\xi_n)$, where $\eta$, $\zeta$, and $\xi$ are the adapted  sequences in \eqref{main-decom}.  Clearly, $y=w^d +w^c +w^r$ and the decomposition is independent of $p$. We provide some general  estimates for the martingales $w^d$, $w^c$, and $w^r$.

\begin{proof}[Norms estimates for the column and row parts] 
 Using the fact that  for every operator $a\in L_2(\M)$ and $n\geq 1$,
$\E_{n-1}|a-\E_{n-1}(a)|^2 = \E_{n-1}|a|^2 -|\E_{n-1}(a)|^2 \leq \E_{n-1}|a|^2$, it is clear that 
\[
\| w^c \|_{\h_p^c} \leq \|\zeta\|_{L_p^{\rm cond}(\M;\ell_2^c)}=\|\sigma_c(\zeta)\|_p
\]
 We claim that for every $1<p<2$,  $\|\sigma_c(\zeta)\|_p^p \leq 1+\kappa_p\|x\|_p^p$ for some suitable constant $\kappa_p$ depending only on $p$.  It suffice to show that for every $1<p<2$ and $N\geq 1$, $\|\sigma_{c,N}(\zeta)\|_p^p \leq 1+\kappa_p\|x_N\|_p^p$.

 Our approach was inspired by  ideas  from  \cite{Jiao-OS-Wu} taking advantage of various properties of projections that are relatives to the Cuculescu projections established  in the previous section. 
Below, we  will use the well-known fact that  for any given operator $a$, its $L_p$-norm can be computed with the integral formula:
\[
\big\|a\big\|_p^p =\int_0^\infty pt^{p-1} \T\big(\ch_{(t,\infty)}(|a|)\big)\ dt.
\]
We begin with  the following elementary estimate:
\begin{equation}\label{estimate-strong}
\big\| \sigma_{c,N}(\zeta)\big\|_p^p \leq 1+ (2^p-1) \sum_{k\geq 0} 2^{pk} \T\big(\ch_{(2^k,\infty)}(\sigma_{c,N}(\zeta))\big).
\end{equation}
According  to  the estimates on the distribution of $\sigma_{c,N}(\zeta)$ stated in Proposition~\ref{distribution-c}, we  have
\begin{equation*}
\begin{split}
\big\| \sigma_{c,N}(\zeta)\big\|_p^p &\leq 1 +  (2^p-1)
 \sum_{k=0}^\infty  2^{pk} \Big(2^{-2k} \big\|e_{k,N }x_N e_{k,N}\big\|_2^2  +  8\ . \ 2^{-k}\T\big(({\bf 1} -e_{k,N})|x_N|\big)\Big)\\
&= 1+  (2^p-1)\sum_{k=0}^\infty  2^{(p-2)k}
\big\| e_{k,N} x_N  e_{k,N}\big\|_2^2  + 8(2^p-1) \sum_{k=0}^\infty  2^{(p-1)k }\T\big( ({\bf 1}-e_{k,N} )|x_N|\big). 
\end{split}
\end{equation*}
The last step is to apply Lemma~\ref{lem-p} and Lemma~\ref{last} to conclude that
\begin{equation}\label{column-p}
\begin{split}
\big\| \sigma_{c,N}(\zeta)\big\|_p^p &\leq 1
  + \Big[\frac{(2^p-1)2^{p^2 +1}}{(1-2^{p-2})(2^{p-1}-1)^p} +8\frac{(2^p-1)2^{(p-1)^2}}{(2^{p-1}-1)^p}\Big] \big\|x_N\big\|_p^p\\
 &=1+  \frac{2^{p^2}}{(2^{p-1}-1)^p}\Big[   \frac{2(2^p-1)}{1-2^{p-2}} + 8 (2^p-1)2^{-2p+1}\Big] \big\|x_N\big\|_p^p\\
 &=1+  \frac{2^{p^2}}{(2^{p-1}-1)^p}\Big[   \frac{2(2^p-1)}{1-2^{p-2}} + 8 (2^{-p+1} -2^{-2p+1})\Big] \big\|x_N\big\|_p^p.
 \end{split}
\end{equation}
As  $2^p-1 \leq 3$ and $2^{-p+1}-2^{-2p+1} \leq 1$, we have
\[
\big\| \sigma_{c,N}(\zeta)\big\|_p^p \leq  1+ \frac{2^{p^2}}{(2^{p-1}-1)^p}\Big[ 8+ \frac{6}{1-2^{p-2}}\Big] \big\|x_N\big\|_p^p.
\]
Taking limits as $N \to \infty$, we conclude that:
\begin{equation}\label{zeta-p}
\big\| w^c \big\|_{\h_p^c}^p \leq  \big\| \sigma_c(\zeta)\big\|_p^p \leq  1+ \frac{2^{p^2}}{(2^{p-1}-1)^p}\Big[ 8+ \frac{6}{1-2^{p-2}}\Big] \big\|x\big\|_p^p.
\end{equation}
Similarly, we may also state
\begin{equation}\label{y^r-p}
\big\| w^r\big\|_{\h_p^r}^p \leq  1+ \frac{2^{p^2}}{(2^{p-1}-1)^p}\Big[ 8+ \frac{6}{1-2^{p-2}}\Big] \big\|x\big\|_p^p.
\end{equation}
\end{proof}

\begin{proof}[Norm estimates  for the diagonal part]

Since conditional expectations are contractions in $L_p(\M)$, we have $\|w^d\|_{\h_p^d} \leq  2\|\eta\|_{L_p(\M \overline{\otimes} \ell_\infty)}$ so we will only work on $\|\eta\|_{L_p(\M \overline{\otimes} \ell_\infty)}$.  As in the column part, it suffices  to estimate the $L_p$-norm of the finite sequence $\eta^{(N)}= (\eta_n)_{1\leq n\leq N}$ for $N\geq 1$.  Elementary calculation shows that 
\[
\big\|\eta^{(N)}\big\|_p^p \leq (2^p-1)  \sum_{k \in \mathbb{Z}} 2^{pk} \Tr\big(\ch_{(2^k,\infty)}(|\eta^{(N)}|)\big).
\]
Next, we divide the right hand side into two parts:
\[
 \sum_{k \in \mathbb{Z}} 2^{pk} \Tr\big(\ch_{(2^k,\infty)}(|\eta^{(N)}|)\big)= \sum_{k \geq 1} 2^{pk} \Tr\big(\ch_{(2^k,\infty)}(|\eta^{(N)}|)\big)+\sum_{-\infty<k\leq 0} 2^{pk} \Tr\big(\ch_{(2^k,\infty)}(|\eta^{(N)}|)\big) =\mathbb{A} + \mathbb{B}.
\]  
From Proposition~\ref{distribution-d}(i), we have
\begin{equation*}
\begin{split}
\mathbb{A} &\leq  4\sum_{k\geq 1} 2^{(p-2)k} \big\| e_{k,N} x_N  e_{k,N}\big\|_2^2 + 28 \sum_{k\geq 1} 2^{(p-1)k} \T\big( ({\bf 1}-e_{k,N} )|x_N|\big) + \sum_{k\geq 1} 2^{pk} \T\big(\ch_{(2^k,\infty)}(|x_1|)\big)\\
&\leq 4\sum_{k\geq 1} 2^{(p-2)k} \big\| e_{k,N} x_N  e_{k,N}\big\|_2^2 + 28 \sum_{k\geq 1} 2^{(p-1)k} \T\big( ({\bf 1}-e_{k,N} )|x_N|\big) + \big\|x_1\big\|_p^p.
\end{split}
\end{equation*}
Using the estimate from Lemma~\ref{last} in the first term,  Lemma~\ref{lem-p} in the second term, and the fact that $\|x_1\|_p\leq \|x_N\|_p$, we arrive at the estimate:
\[
\mathbb{A} \leq \Big[ 4 \frac{2^{p^2 +1}}{(1-2^{p-2})(2^{p-1}-1)^p}+  28\frac{2^{(p-1)^2}}{(2^{p-1}-1)^p} + 1 \Big] \big\|x_N\big\|_p^p.
\]
We may  record  after simplification that:
\begin{equation}\label{positive}
(2^p-1)\mathbb{A} \leq 
\frac{ 2^{p^2}}{(2^{p-1}-1)^p}\Big[15(2^p-1) +  \frac{8(2^p-1)}{1-2^{p-2}}  \Big] \big\|x_N\big\|_p^p.
\end{equation}
Now we estimate  $\mathbb{B}$. Since the trace $\T$ is normalized, it follows from Proposition~\ref{distribution-d}(ii)   that for every $k\leq 0$,
\begin{align*}
\Tr\big(\ch_{(2^k,\infty)}(|\eta^{(N)}|)\big) &\leq 1+\sum_{i \geq 0} \T({\bf 1} -e_{i,N})\\
&\leq 1+ \sum_{i\geq 0} \sum_{j\geq i} \T({\bf 1} -q_N^{(2^j)}).
\end{align*}
We claim that $\T({\bf 1}-q_N^{(2^j)}) \leq  2^{-pj} \|x_N\|_p^p$. This follows from the following fact (whose verification can be found in the proof  of \cite[Proposition~1.4]{PR}): for $l\geq 1$,
\[
2^{j}(q_{l-1}^{(2^j)} -q_{l}^{(2^j)})\leq |(q_{l-1}^{(2^j)} -q_{l}^{(2^j)})(q_{l-1}^{(2^j)} x_l q_{l-1}^{(2^j)})(q_{l-1}^{(2^j)} -q_{l}^{(2^j)})|.
\]
This implies that  for every $1\leq l \leq N$, 
\begin{align*}
\T(q_{l-1}^{(2^j)} -q_{l}^{(2^j)}) &\leq 2^{-jp} \big\| (q_{l-1}^{(2^j)} -q_{l}^{(2^j)})x_l(q_{l-1}^{(2^j)} -q_{l}^{(2^j)})\big\|_p^p \\
&\leq  2^{-jp} \big\| (q_{l-1}^{(2^j)} -q_{l}^{(2^j)}) x_N(q_{l-1}^{(2^j)} -q_{l}^{(2^j)})\big\|_p^p.
\end{align*}
Taking the summation over $l$ gives the claim. We now have for  every $k\leq 0$,
\begin{align*}
\Tr\big(\ch_{(2^k,\infty)}(|\eta^{(N)}|) \big) &\leq  1+ \sum_{i\geq 0} \sum_{j\geq i} 2^{-jp} \big\|x_N\big\|_p^p\\
&=1+ \Big[ \frac{2^{2p}}{(2^p-1)^2}\Big] \big\|x_N \big\|_p^p.
\end{align*}
%%%%%%
%%%%%%%%
%%%%%%%%%%%
 It then follows that 
\begin{equation}\label{negative}
(2^p-1)\mathbb{B} \leq  2^p + \Big[ \frac{2^{3p}}{(2^p-1)^2}\Big] \big\|x_N \big\|_p^p.
\end{equation}
Combining the two estimates \eqref{positive} and \eqref{negative}, and taking limits as $N\to \infty$, we arrive at
\[
\big\| \eta\big\|_p^p \leq   2^p + 
\frac{2^{p^2}}{(2^{p-1}-1)^p}\Big[15(2^p-1) +  \frac{8(2^p-1)}{1-2^{p-2}} +\frac{2^{3p}  (2^{p-1}-1)^p}{2^{p^2}(2^p-1)^2}\Big] \big\| x\big\|_p^p.
\]
Using $2(2^{p-1}-1)^p/ 2^{p^2} \leq 1$ and $1\leq 2^p -1 \leq 3$, we have after simplification of the right hand side  that:
\begin{equation}\label{eta-p}
\big\|w^d\big\|_{\h_p^d}^p\leq(2\big\| \eta\big\|_p)^p \leq 2^{2p} + \frac{2^{p^2+p}}{(2^{p-1}-1)^p}\Big[77 +  \frac{24}{1-2^{p-2}}  \Big] \big\|x\big\|_p^p.
\end{equation}
Thus we obtain an estimate for the diagonal part.

%%%%%%%%
%%%%%%%%

\medskip

{\bf Step~2.} We should emphasize  that the estimates in \eqref{zeta-p}, \eqref{y^r-p}, and \eqref{eta-p} in the previous step are non-homogeneous. In this step, we will remedy this imperfection.

For each $m\in \mathbb{N}$, consider the martingale  $x^{(m)}= m x$ and $y^{(m)} =my$. Then $x^{(m)}$ is a self-adjoint  $L_2$-martingale and $y^{(m)}$ is weakly differentially subordinate to $x^{(m)}$.  According to Step~1, $y^{(m)}$ admits a decomposition $y^{(m)} =y^{(m,d)} +y^{(m,c)} + y^{(m,r)}$ satisfying the norm estimates for every $1<p<2$:
\begin{equation}\label{m-p}
\begin{split}
\big\| \frac{y^{(m,c)}}{m}\big\|_{\h_p^c}^p &\leq  \frac{1}{m^p}+ \frac{2^{p^2}}{(2^{p-1}-1)^p}\Big[ 8+ \frac{6}{1-2^{p-2}}\Big] \big\|x\big\|_p^p;\\
\big\| \frac{y^{(m,r)}}{m}\big\|_{\h_p^r}^p &\leq  \frac{1}{m^p}+ \frac{2^{p^2}}{(2^{p-1}-1)^p}\Big[ 8+ \frac{6}{1-2^{p-2}}\Big] \big\|x\big\|_p^p;\\
\big\|\frac{y^{(m,d)}}{m}\big\|_{\h_p^d}^p &\leq \frac{2^{2p}}{m^p} + \frac{2^{p^2+p}}{(2^{p-1}-1)^p}\Big[77 +  \frac{24}{1-2^{p-2}}  \Big] \big\|x\big\|_p^p.
\end{split}
\end{equation}
Consider the sequence $\displaystyle{\big\{{y^{(m,c)}}/{m}\big\}_{m\geq 1}}$ in $\h_1^c(\M)$. Since it is bounded in $\h_p^c(\M)$ for all $1<p<2$, by taking a subsequence if necessary, we may assume that it admits a weak-limit  in $\h_1^c(\M)$. Similarly, we may assume that $\displaystyle{\big\{{y^{(m,r)}}/{m}\big\}_{m\geq 1}}$ admits a weak-limit in $\h_1^r(\M)$. 
For the diagonal part, consider the sequence $\displaystyle{\big\{{dy^{(m,d)}}/{m}\big\}_{m\geq 1}}$ in $\ell_1(L_1(\M)) \subset \ell_2(L_1(\M))$. Since for any $1<p<2$,
$\ell_p(L_p(\M)) \subset \ell_2(L_1(\M))$ 
and  $\displaystyle{\big\{{dy^{(m,d)}}/{m}\big\}_{m\geq 1}}$ is bounded in $\ell_p(L_p(\M))$, identical argument allows  us to also assume that  $\displaystyle{\big\{dy^{(m,d)}/m\big\}_{m\geq 1}}$ converges weakly in $\ell_2(L_1(\M))$. Set
$y^c:= w-\lim_{m\to \infty} {y^{(m,c)}}/{m}$ in  $\h_1^c(\M)$, $y^r:= w-\lim_{m\to \infty} {y^{(m,r)}}/{m}$ in  $\h_1^r(\M)$, and $dy^d:= w-\lim_{m\to \infty} { dy^{(m,d)}}/{m}$ in  $\ell_2(L_1(\M))$. Then, we clearly have $y=y^c + y^r +y^d$ (as  sequences in $L_1(\M)$) and this decomposition is independent of $p$. 

Now fix $1<p<2$. There exists an increasing sequence of integers $(m_j)_{j\geq 1}$ (which may depend on $p$) so  that 
for  $\dagger \in \{c,r,d\}$,  $\displaystyle{\big\{{y^{(m_j,\dagger)}}/{m_j}\big\}_{j\geq 1}}$ converges weakly in $\h_p^{\dagger}(\M)$. Since  the inclusions $\h_p^c(\M) \subset \h_1^c(\M)$, $\h_p^r(\M) \subset \h_1^r(\M)$, and $\ell_p(L_p(\M)) \subset \ell_2(L_1(\M))$  are  continuous, such weak-limit must be equal to $y^\dagger$.
It then follows from \eqref{m-p} that
\begin{equation*}
\begin{split}
\big\| y^c\big\|_{\h_p^c}^p &\leq   \frac{2^{p^2}}{(2^{p-1}-1)^p}\Big[ 8+ \frac{6}{1-2^{p-2}}\Big] \big\|x\big\|_p^p;\\
\big\| y^r\big\|_{\h_p^r}^p &\leq   \frac{2^{p^2}}{(2^{p-1}-1)^p}\Big[ 8+ \frac{6}{1-2^{p-2}}\Big] \big\|x\big\|_p^p;\\
\big\|y^d\big\|_{\h_p^d}^p &\leq  \frac{2^{p^2+p}}{(2^{p-1}-1)^p}\Big[77 +  \frac{24}{1-2^{p-2}}  \Big] \big\|x\big\|_p^p.
\end{split}
\end{equation*}

We can now conclude from combining these three inequalities   that 
\[
\|y^d \|_{\h_p^d} +\|y^c \|_{\h_p^c} + \|y^r \|_{\h_p^r}   \leq c_p \big\|x \big\|_p
\]
where $c_p$  is the constant in the statement of the theorem.  This completes the proof.
\end{proof}

The construction used in the proof of the preceding theorem also provides estimates of the $L_2$-norms. The following corollary improves \cite[Theorem~6.1.8]{Junge-Perrin} on three fronts: first, it is valid for    
 weakly differentially subordinate martingales, second the decomposition is valid simultaneously for all  $1<p \leq 2$, and finally  it gives a much simpler estimate on the $L_2$-norms. It also improves the simultaneous decomposition from \cite[Theorem~2.10]{RW}.  We should note that the approaches used in  both \cite{Junge-Perrin} and \cite{RW} are interpolations as opposed to the constructive approach   used here.

\begin{corollary}\label{strong-type-cor} Suppose  $x$ is a self-adjoint $L_2$-bounded martingale  and $y$ is a self-adjoint martingale that is weakly differentially subordinate to $x$. Then there exists three martingales $y^d$, $y^c$, and $y^r$  (depending only on $x$ and $y$) such that:
\begin{enumerate}[{\rm (i)}]
\item $y= y^d +y^c +y^r$;
\item for every $1<p<2$, $\big\|y^d\big\|_{\h_p^d} +\big\|y^c\big\|_{\h_p^c} +\big\|y^r\big\|_{\h_p^r}\leq C_p \big\|x \big\|_p$;
\item $\max\big\{\big\|y^d \big\|_2, \big\|y^c \big\|_2, \big\|y^r \big\|_2\big\} \leq 5\big\|x\big\|_2$.
\end{enumerate}
Here  $C_p \leq C(p-1)^{-1} $ as $p\to 1$.
\end{corollary}
\begin{proof} We use the decomposition in Theorem~\ref{strong-type-1} to get $(i)$ and $(ii)$. It is clear from \eqref{main-decom} that $\max\{\|y^c\|_2, \|y^r\|_2\} \leq 2 \|y\|_2 \leq  2\|x\|_2$. Triangle inequality then gives $\|y^d\|_2\leq 5\|y\|_2 \leq  5\|x\|_2$.
\end{proof}
%%%%%%%%%%%%%%%
%%%%%%%%%%%%%%

We now present the corresponding result for Hardy spaces norms.

\begin{theorem}\label{strong-type-2}  
Suppose that $x$ is a  self-adjoint $L_2$-martingale and $y$ is a  self-adjoint martingale that  is weakly  differentially subordinate to $x$. Then
 there exist
two martingales $y^c$ and $y^r$  (depending only on  $x$ and $y$) such
that
$y= y^c +y^r$ and for every $1<p<2$,
\[
 \big\| y^c\big\|_{\H_p^c}  +\big\| y^r\big\|_{\H_p^r} \leq c_p \big\| x \big\|_p
\]
with $c_p= O((p-1)^{-1})$ when $p \to 1$.
\end{theorem}

\begin{proof} We assume without loss of generality that $x$ is a finite martingale and   consider the decomposition $y=y^c + y^r$ provided by \eqref{equation-square}. It is enough to consider just the column part since the row part can be proved using  identical argument. To estimate $\|y^c\|_{\H_p^c}$, we use as before the elementary inequality:
\[
\big\|y^c\big\|_{\H_p^c}^p \leq 1+  (2^p-1)\sum_{k\geq 0} 2^{pk} \T\Big(\ch_{(2^k, \infty)}\big(S_c(y^c)\big) \Big).
\]
Using estimates on  the distribution function of $S_c(y^c)$ from Remark~\ref{distribution-similar}, the proof is now just a notational adjustment of the proof for $\|\sigma_c(\zeta)\|_p$ in Theorem~\ref{strong-type-1}. Details are left to the reader.
\end{proof}

%%%%%%%%%%%%%%%
%%%%%%%%%%%%%%%%
%%%%%%%%%%%%%%%%%
%%%%%%%%%%%%%%%%%

\subsection{Comparisons of norms}
In this subsection, 
we will  compare the three strong type results  Theorem~\ref{strong-type-1}, 
Theorem~\ref{strong-type-2}, and \cite[Theorem~5.1(i)]{Jiao-OS-Wu}. This  will be achieved by comparing for $1\leq p<2$, the $L_p$-norms, the conditioned Hardy space norms, and the Hardy space norms  in the style of noncommutative Burkholder inequalities and noncommutative Burkholder-Gundy inequalities (\cite{JX,PX}).

We take the opportunity  to present  our results with their respective sharp constants which are interesting questions on their own. 
The next  theorem contains 
noncommutative analogues of   results of Wang in \cite{Wang} (see also \cite[Theorem~8.19]{OS}). These inequalities are  new  and could be of independent interest.  We should note here that except for the constants,   the range  $1 \leq p<2$ are  known as they are parts of the noncommutative Burkholder inequalities and the noncommutative Burkholder-inequalities (see for instance, \cite{JX2}). Our approach below is very different from Wang's proof. It
was inspired by an argument used in \cite{Bekjan-Chen-Perrin-Y}  to  describe  an equivalent  quasi-norm on $\h_p^c(\M)$ when $0<p<2$ which in turn was adapted from an argument  due to Herz \cite{Herz} for the classical case. We now state the main result of this subsection: 
\begin{theorem}\label{L-h}
Let  $0<p\leq 2$. For every $x \in \h_p^c(\M)$, the following two inequalities hold:
\[
\big\|x \big\|_p \leq  \sqrt{2/p} \big\| x \big\|_{\h_p^c} 
\]
and
\[
\big\|x \big\|_{\H_p^c} \leq  \sqrt{2/p} \big\| x \big\|_{\h_p^c}. 
\]
In each of the two inequalities, the constant $\sqrt{2/p}$ is the best possible. 
\end{theorem}
For the proof we need the following lemma:
\begin{lemma}[{\cite[Lemma~3.1]{Bekjan-Chen-Perrin-Y}}]\label{integral}
Let $f $ be a function in $C^1(\mathbb{R}_+)$ and  $a, b \in \M_+$, then
\[
\T\big(f(a+b)-f(a)\big) =\T\left( \int_0^1 f'(a+tb)b\ dt \right).
\]
\end{lemma}

\begin{proof}[Proof of Theorem~\ref{L-h}] Let  $W$ be the collection of all sequences of positive operators $(w_n)_{n\geq 1}$ such that $\{w_n^{-1+\frac{2}{p}}\}_{n\geq 1}$  is nondecreasing with each $w_n \in L_1(\M_{n-1})$,  is invertible with bounded inverse, and satisfies $\|w_n\|_1\leq 1$. For an $L_2$-martingale $x$, we set 
\[
N_p^c(x)=\inf\left\{\Big[ \T\big( \sum_{n\geq 1} w_n^{1-\frac{2}{p}} |dx_n|^2\big) \Big]^{1/2}: (w_n)_{n\geq 1} \in W \right\}.
\]
We estimate $N_p^c(x)$ from above and from below. The estimate from above is already included in the proof of \cite[Proposition~3.2]{Bekjan-Chen-Perrin-Y} but we include the argument for completeness.  Let $x \in L_2(\M)$ with  $\|x\|_{\h_p^c}=1$. By approximation, we may assume that $x \in \M$ and  $s_{c,n}(x)$ is invertible with bounded inverse for every $n\geq 1$.

For $n\geq 1$, take $w_n= s_{c,n}^p(x)$. Then,  $(w_n)_{n\geq 1} \in W$. Since the sequence $(w_n)_{n\geq 1}$ is predictable, we have
\begin{align*}
\T\big( \sum_{n\geq 1} w_n^{1-\frac{2}{p}} |dx_n|^2\big) &=\T\big( \sum_{n\geq 1} w_n^{1-\frac{2}{p}}\E_{n-1}( |dx_n|^2)\big)\\
&=\T\big( \sum_{n\geq 1} s_{c,n}^{p-2}(x) (s_{c,n}^2(x)- s_{c,n-1}^2(x))\big).
\end{align*}
Applying Lemma~\ref{integral} with  $f(t)=t^{p/2}$, $a+b=s_{c,n}^2(x)$ and $a=s_{c,n-1}^2(x)$, we obtain
\begin{align*}
\frac{p}{2}&\T\big( s_{c,n}^{p-2}(x) \big[s_{c,n}^2(x)- s_{c,n-1}^2(x)\big]\big)\\ &\leq \T\left(\int_0^1 \frac{p}{2}\big[ s_{c,n-1}^2(x) + t \big(s_{c,n}^2(x)-s_{c,n-1}^2(x) \big)\big]^{\frac{p}{2}-1} \big[s_{c,n}^2(x)-s_{c,n-1}^2(x) \big] \ dt \right)\\ 
&= \T\big(s_{c,n}^p(x)-s_{c,n-1}^p(x)\big)
\end{align*}
where we have used the fact that the function $t \mapsto t^{\frac{p}{2}-1}$ is operator nonincreasing when $0<p<2$. Taking summation over $n$ leads to 
\begin{equation}\label{above}
N_p^c(x)^2 \leq \frac{2}{p} \T\big(s_c^p(x)\big) \leq \frac{2}{p}.
\end{equation}

Now we estimate $N_p^c(x)$ from below. Fix $(w_n)_{n\geq 1} \in W$.  Since $(w_n)_{n\geq 1}$ is a predictable sequence, we have
\begin{align*}
\T\big(  w_n^{1-\frac{2}{p}} |dx_n|^2\big) &=\T\big(  w_n^{1-\frac{2}{p}} \E_{n-1}(|dx_n|^2)\big)\\
&=\T\big(  w_n^{1-\frac{2}{p}} \E_{n-1}(|x_n -x_{n-1}|^2)\big)\\
&=\T\big(  w_n^{1-\frac{2}{p}} \E_{n-1}(|x_n|^2 -|x_{n-1}|^2)\big)\\
&=\T\big(  w_n^{1-\frac{2}{p}} (|x_n|^2 -|x_{n-1}|^2)\big)\\
&= \T\big(  w_n^{1-\frac{2}{p}} |x_n|^2 \big)- \T\big(w_n^{1-\frac{2}{p}}|x_{n-1}|^2\big).
\end{align*}
Since $(w_n^{1-\frac{2}{p}})_{n\geq 1}$ is decreasing, we have $ \T\big(  w_n^{1-\frac{2}{p}} |x_n|^2 \big)\geq  \T\big(  w_{n+1}^{1-\frac{2}{p}} |x_n|^2 \big)$. This implies that
\[
\T\big(  w_n^{1-\frac{2}{p}} |dx_n|^2\big) \geq \T\big(  w_{n+1}^{1-\frac{2}{p}} |x_n|^2 \big)- \T\big(w_n^{1-\frac{2}{p}}|x_{n-1}|^2\big).
\]
Taking summation over $n$ (with $x_0=0$), we have for every $k\geq 1$,
\[
\T\big(  w_{k+1}^{1-\frac{2}{p}} |x_k|^2\big)\leq   \T\big( \sum_{n\geq 1} w_n^{1-\frac{2}{p}} |dx_n|^2\big). 
\]
Using H\"older's inequality, we get 
\begin{align*}
\big\|x_k\big\|_p &\leq  \big\|w_{k+1}\big\|_1^{-\frac{1}{2}+\frac{1}{p}} \big\| x_k w_{k+1}^{\frac{1}{2}-\frac{1}{p}}\big\|_2\\
&\leq \big[\T\big( \sum_{n\geq 1} w_n^{1-\frac{2}{p}} |dx_n|^2\big)\big]^{1/2}.
\end{align*}
Taking the infimum over all $(w_n)_{n\geq 1} \in W$ and the supremum over all $k\geq 1$, we obtain 
\begin{equation}\label{below}
\big\| x\big\|_p \leq N_p^c(x). 
\end{equation}
The first inequality  in Theorem~\ref{L-h} clearly follows from combining \eqref{above} and \eqref{below}.

The second inequality can be deduced using similar  argument. Indeed,  for every $n\geq 1$, 
\begin{align*}
\T\big(  w_n^{1-\frac{2}{p}} |dx_n|^2\big) &=\T\big(  w_n^{1-\frac{2}{p}} [S_{c,n}^2(x)- S_{c,n-1}^2(x)]\big)\\
&\geq \T\big(  w_{n+1}^{1-\frac{2}{p}} S_{c,n}^2(x)\big)- \T\big( w_{n}^{1-\frac{2}{p}} S_{c,n-1}^2(x)\big).
\end{align*}
The rest of the proof is identical with $S_{c,n}^2(x)$ in place of $|x_n|^2$. The fact that the constant $\sqrt{2/p}$ is sharp is already the case for the commutative martingales as shown in \cite[Theorem~1]{Wang}.
\end{proof}

We  can also make the comparisons at the level  of   mixed Hardy spaces and  mixed conditioned Hardy spaces. The following follows directly  from Theorem~\ref{L-h} and its row version.  
\begin{corollary}\label{H-h}
Let  $0<p \leq 2$. For every $x \in \h_p(\M)$,  we have
\[
\big\|x \big\|_{p} \leq  c_p \big\| x \big\|_{\h_p}
\]
where
$c_p=\begin{cases} 3^{\frac{1-p}{p}}\sqrt{2/p} &{\text if} \ 0<p<1\\
\sqrt{2/p} &{\text if} \ 1\leq p \leq 2.
\end{cases}
$

Similarly,  we have 
\[
\big\|x \big\|_{\H_p} \leq  c_p' \big\| x \big\|_{\h_p}
\]
where
$c_p'=\begin{cases} 2^{\frac{1-p}{p}}\sqrt{2/p} &{\text if} \ 0<p<1\\
\sqrt{2/p} &{\text if} \ 1\leq p \leq 2.
\end{cases}
$
\end{corollary}

\begin{remark}\label{L-H}
From the noncommutative Davis decomposition (\cite{Junge-Perrin, Perrin,Ran-Wu-Xu}), we may also deduce from  Theorem~\ref{L-h} that for $1\leq p< 2$ and $x \in \H_p(\M)$, we have
\begin{equation}\label{e:LH}
\big\| x \big\|_p \leq c_p \big\| x\big\|_{\H_p}
\end{equation}
with $c_p=\big( (\delta_{p/(2-p)}')^{1/2} +1\big)\sqrt{2/p} +2/p$ where $\delta_{p/(2-p)}'$ denotes the constant from the noncommutative dual Doob inequality relative  to the index $p/(2-p)$ (\cite{Ju,JX2}).
This should be compared with \cite[Corollary~4.3]{JX}. We should point out that even for classical martingales, \eqref{e:LH}  is not valid for $0<p<1$ (see \cite[Example~8.1]{BG}).
\end{remark}
 We refer to \cite{Ran-Wu-Xu} for the definition of the space $\h_p^{1_c}(\M)$ used below.  In order to see \eqref{e:LH}, let  $\epsilon>0$ and  fix a decomposition $x=x^{(1)} + x^{(2)}$ so that
\[
\big\| x^{(1)}\big\|_{\H_p^{c}} +\big\| x^{(2)}\big\|_{\H_p^{r}} \leq \big\| x\big\|_{\H_p} +\epsilon.
\]
Inspecting the proof of 
 the noncommutative Davis decomposition  from \cite{Junge-Perrin}, we have two martingales $x^{(1,d)}$ and $x^{(1,c)}$ with $x^{(1)} =x^{(1,d)} +x^{(1,c)}$, 
$\|x^{(1,d)}\|_{\h_p^{1_c}} \leq  \big( (\delta_{p/(2-p)}')^{1/2} +1\big)\sqrt{2/p} \|x^{(1)}\|_{\H_p^c}$  and $\|x^{(1,c)}\|_{\h_p^{c}} \leq  \sqrt{2/p} \|x^{(1)}\|_{\H_p^c}$. For more details on  how the above constants are deduced, we refer to the argument in the proof of \cite[Lemma~6.1.4]{Junge-Perrin}. 

 It is clear that $\|x^{(1,d)}\|_{p}  \leq \|x^{(1,d)}\|_{\h_p^{1_c}} \leq  \big( (\delta_{p/(2-p)}')^{1/2} +1\big)\sqrt{2/p}\|x^{(1)}\|_{\H_p^c}$. By Theorem~\ref{L-h}, we also have 
$\|x^{(1,c)}\|_{p} \leq \sqrt{2/p}\|x^{(1,c)}\|_{\h_p^{c}} \leq  (2/p) \|x^{(1)}\|_{\H_p^c}$.
Combining the two estimates, we clearly get:
\[
\| x^{(1)} \|_p \leq \big[\big( (\delta_{p/(2-p)}')^{1/2} +1\big)\sqrt{2/p} + 2/p\big] \|x^{(1)}\|_{\H_p^c}.
\]
Similar estimate applies to $x^{(2)}$.
\qed

\medskip

Using Corollary~\ref{H-h}, we see that Theorem~\ref{strong-type-1} implies both \cite[Theorem~5.1(i)]{Jiao-OS-Wu} and Theorem~\ref{strong-type-2} while still maintaining the optimal  orders of the constants when $p\to 1$. It is worth   pointing out  that using the noncommutative Burkholder/Rosenthal inequalities, one can still deduce Theorem~\ref{strong-type-1} from \cite[Theorem~5.1(i)]{Jiao-OS-Wu}. However, we no longer able to maintain the correct order of constants if we follow  that route.  Therefore, we  may treat Theorem~\ref{strong-type-1} as  the strongest  of the three strong-type results. On the other hand, the decomposition used in the proof of Theorem~\ref{strong-type-1}  is far more complex than the one used in Theorem~\ref{strong-type-2} and \cite[Theorem~5.1(i)]{Jiao-OS-Wu} does not require any decomposition so it is still justified to have separate proofs for Theorem~\ref{strong-type-2} and \cite[Theorem~5.1(i)]{Jiao-OS-Wu}.

\medskip

At the time of this writing, we do not know if one  can compare the corresponding weak-type results  (Theorem~\ref{main-weak}, Theorem~\ref{main-weak-S}, and \cite[Theorem~4.1]{Jiao-OS-Wu}). For instance, unlike the case of $0<p \leq 2$ stated in Theorem~\ref{L-h}, we do not know if there exists an absolute  constant $C$ with $\|x\|_{1,\infty} \leq C\|s_c(x)\|_{1,\infty}$  (or $\|S_c(x)\|_{1,\infty} \leq C\|s_c(x)\|_{1,\infty}$) for every martingale $x$. Another obstacle for potential comparisons is the fact that it is still unknown if the decomposition in Theorem~\ref{main-weak} can be achieved with martingale difference sequences as opposed to just adapted sequences.
As a result, it is  unclear at this point if 
Theorem~\ref{main-weak} implies \cite[Theorem~4.1]{Jiao-OS-Wu} or Theorem~\ref{main-weak-S}.

\bigskip
\begin{remark}
It is important to observe that  the argument used  to prove Theorem~\ref{strong-type-1} (if one only considers the special case  $y=x$) and Corollary~\ref{H-h}  provide direct and constructive proofs of the noncommutative Burkholder/Rosenthal inequalities for $1<p<2$ with  optimal order of the constants when $p\to 1$. To the best of our knowledge, all available proofs prior to this point were either by duality  or relaying on interpolation techniques. Likewise, Theorem~\ref{strong-type-2} and Corollary~\ref{H-h} provide direct and constructive proofs of the noncommutative Burkholder-Gundy inequalities for $1<p<2$  with  optimal order of the constants when $p\to 1$.
\end{remark}
%%%%%%%%%%%%%%%

\section{Remarks and open problems}

\subsection{Modular inequalities}

In this subsection, we  will demonstrate  that our method of proof in   Subsection~\ref{weak-type} is general enough  to cover some special  cases  of  noncommutative moment inequalities associated with Orlicz functions.   By an  \emph{ Orlicz function}  $\Phi$ on $[0,\infty)$, we mean a continuous, increasing, and convex function such that $\Phi(0)=0$ and $\lim_{t\to \infty} \Phi(t)=\infty$. We will assume throughout that 
 $\Phi$ is an Orlicz function satisfying the $\Delta_2$-condition, that is, for some constant $C>0$,
\begin{equation}
\Phi(2t) \leq C \Phi(t), \quad  t\geq 0.
\end{equation}
Given $1\leq p\leq q < \infty$,  we recall that an Orlicz function $\Phi$  is said to be \emph{$p$-convex} if the function $t \mapsto\Phi(t^{1/p})$ is convex, and to be  \emph{$q$-concave} if the function $t\mapsto \Phi(t^{1/q})$ is concave. The function $\Phi$ satisfies the $\Delta_2$-condition if and only if  it is $q$-concave for some $q<\infty$. If $\Phi$ is $q$-concave for $1<q<\infty$ then the function $t \mapsto t^{-1}\Phi(t)$ is increasing and the function $t\mapsto t^{-q}\Phi(t)$ is decreasing.
We denote by $L_\Phi$ the associated Orlicz function space. 
 We refer to \cite{Kras-Rutickii, Maligranda2}   for backgrounds on  Orlicz functions and  Orlicz spaces.

Given an operator $x \in L_0(\M,\T)$ and an Orlicz function $\Phi$, we may define $\Phi(|x|)$ using functional calculus. That is, if $|x|=\int_0^\infty s \ de_s^{|x|}$ is its spectral decomposition, then 
\[
\Phi(|x|) =\int_0^\infty \Phi(s)\ de_s^{|x|}.
\]
In this case, $\Phi(|x|) \in L_0(\M,\T)$ and its trace  $\T( \Phi(|x|))$  is referred to as the \emph{$\Phi$-moment} of  the operator $|x|$. 
For more information  and background  on moment inequalities involving convex functions, we refer to  \cite{Bekjan-Chen, Bekjan-Chen-Ose, Dirksen-Ricard, Jiao-Sukochev-Zanin,Jiao-Sukochev-Zanin-Zhou,RW2, Ran-Wu-Xu}.

 Below, the following two basic facts will be used. First, the Orlicz function $\Phi$  has a representation
\[
\Phi(t)=\int_0^t \varphi(s)\ ds,\quad t>0
\]
where $\varphi$ is  the right derivative of $\Phi$. Second, if $a$ is a positive $\T$-measurable operator then according to \cite[Corollary~2.8]{FK}, we have:
\[
\T\big(\Phi(a)\big)=\int_0^\infty \Phi(\mu_t(a))\ dt.
\]

As noted earlier, we aim to provide  modular extensions of the weak-type $(1,1)$  results stated in Theorem~\ref{main-weak} and Theorem~\ref{main-weak-S}. We begin with the conditioned version which reads as follows:

\begin{theorem}\label{open-B}
Let $\Phi$ be an Orlicz function that is $q$-concave for some $1\leq q<2$. There exists a constant $C_\Phi$ so that if
$x$ is a  self-adjoint $L_2$-martingale  that is bounded in $L_\Phi(\M)$ and $y$ is a  self-adjoint martingale that  is weakly  differentially subordinate to $x$, then
 there exist
three adapted sequences $\eta=(\eta_n)_{n\geq 1}$, $\zeta=(\zeta_n)_{n\geq 1}$, and $\xi=(\xi_n)_{n\geq 1}$ such
that
$dy=\eta + \zeta +\xi$ and satisfy for every $\lambda>0$:
\[
 \Phi(\lambda)\Tr\big( \ch_{(\lambda, \infty)}(|\eta|) \big) + \Phi(\lambda)\T\Big( \ch_{(\lambda,\infty)}\big(\sigma_c(\zeta)\big) \Big) + \Phi(\lambda)\T\Big( \ch_{(\lambda,\infty)}\big(\sigma_r(\xi)\big) \Big)
\leq C_\Phi  \sup_{n\geq 1}\T\big(\Phi(|x_n|)\big).
\]
\end{theorem}

The modular extension of Theorem~\ref{main-weak-S} is formulated as follows:

\begin{theorem}\label{open-BG}
Let $\Phi$ be an Orlicz function that is $q$-concave for some $1\leq q<2$. There exists a constant $C_\Phi$  such that if $x$ is a  self-adjoint $L_2$-martingale that is bounded in $L_\Phi(\M)$ and $y$ is a self-adjoint  martingale that is weakly differentially subordinate  to $x$, then there exist two martingales $y^c$ and $y^r$ such that $y= y^c + y^r$ and for every $\lambda>0$:
\[
\Phi(\lambda) \T\Big( \ch_{(\lambda,\infty)}\big(S_c(y^c)\big) \Big) + \Phi(\lambda)\T\Big( \ch_{(\lambda,\infty)}\big( S_r(y^r)\big) \Big)
\leq C_\Phi  \sup_{n\geq 1}\T\big(\Phi(|x_n|)\big).
\]
\end{theorem}

\begin{remark}
The preceding  two theorems are  new even for the case $y=x$. In that particular situation,
the only relevant case  is when $\Phi$ is such that $\sup\{ p: \Phi\ \text{is $p$-convex}\}=1$. Since when $\Phi$ is $p$-convex  and $q$-concave for $1<p\leq q<\infty$, stronger results  on $\Phi$-moment Burkholder/Rosenthal  inequalities  (\cite{RW2,Ran-Wu-Xu}) and $\Phi$-moment Burkholder-Gundy inequalities (\cite{Bekjan-Chen, Dirksen-Ricard, Jiao-Sukochev-Zanin-Zhou}) are  available.

Theorem~\ref{open-B} and   Theorem~\ref{open-BG}   partially answer  open problems from \cite[Remark~6.1]{Bekjan-Chen-2}. We do not know if the assumption $1\leq q<2$ can be removed.
\end{remark}

For the proofs,  we need two auxiliary lemmas about the Cuculescu projections when the associated martingale is bounded in a noncommutative  Orlicz space.  The next lemma may be viewed as a $\Phi$-moment version of Proposition~\ref{Cuculescu}(iv). We use this as a  convex function companion of Lemma~\ref{lem-p}. 
\begin{lemma}\label{Phi-1}
Let $\Phi$ be an Orlicz  function  that satisfies the $\Delta_2$-condition.  If $x=(x_n)_{n\geq 1}$  is a self-adjoint martingale that is bounded in $L_\Phi(\M)$,  then for every $N\geq 1$ and $k\geq0$,
\[
 \Phi(2^k)2^{-k}\tau\big(({\bf 1}-e_{k,N})|x_N|\big) \leq c_\Phi \T\big(\Phi(|x_N|)\big).
\]
\end{lemma}
\begin{proof} Under the $\Delta_2$-condition, $\Phi$ is $q$-concave for some $q<\infty$. If $q=1$, then there is nothing to prove so we will assume that $\inf\{q: \Phi \ \text{is $q$-concave}\}>1$. We note  from the representation of $\Phi$ and the $q$-concavity that for every $u\geq 0$,
\begin{equation}\label{fact-phi}
\Phi(u) \leq u\varphi(u) \leq  q \Phi(u)
\end{equation}
where $\varphi$ denotes  the right derivative of $\Phi$.  Next, we observe that if $\a \in \R_+$ and $\pi \in \M$ is a projection, then 
\[
\Phi(\a \pi)  = \Phi(\a)\pi.
\]
We will make use of the classical Young inequality (see for instance, \cite[Chap~I]{ Maligranda2}) which states that for every $u,v \in \mathbb{R_+}$,
\[
uv \leq \Phi(u) + \Phi^*(v).
\]
  We are now ready to present the proof. Since $q >1$,  the complementary function $\Phi^*$ is a $p$-convex  Orlicz function for some   $1<p<\infty$. We may choose
 $t_\Phi$ so that for every $u \in \R_+$, 
\begin{equation}\label{phi-small}
\Phi^*(t_\Phi u) \leq (2q)^{-1}\Phi^*(u).
\end{equation}
First, we have from \eqref{fact-phi} that
\[
\Phi(2^k)2^{-k}\T\big(({\bf 1}-e_{k,N})|x_N|\big)\leq  \phi(2^k)\T\big(({\bf 1}-e_{k,N})|x_N|\big).
\]
Next,  we have from \cite[Theorem~4.2(iii)]{FK} that:
\begin{equation*}
\begin{split}
\phi(2^k)\T\big(({\bf 1}-e_{k,N})|x_N|\big)&= \int_0^1 \mu_t\big(\phi(2^k)({\bf 1}-e_{k,N})|x_N|\big)\ dt\\
&\leq \int_0^1 \mu_t\big( t_\Phi \phi(2^k)({\bf 1}-e_{k,N})\big)\  . \ \mu_t\big(t_\Phi^{-1} |x_N|\big) \  dt.
\end{split} 
\end{equation*}
 Applying  the Young inequality, we further get
\begin{equation*}
\begin{split}
\Phi(2^k)2^{-k}\T\big(({\bf 1}-e_{k,N})|x_N|\big) &\leq \int_0^1  \Phi^*\big[ \mu_t\big( t_\Phi \phi(2^k)({\bf 1}-e_{k,N})\big)\big] +\Phi\big[\mu_t\big(t_\Phi^{-1} |x_N|\big)\big] \ dt \\
&=\int_0^1  \Phi^*\big[ \mu_t\big( t_\Phi \phi(2^k)({\bf 1}-e_{k,N})\big)\big] \ dt + \int_0^1 \Phi\big[\mu_t\big(t_\Phi^{-1} |x_N|\big)\big] \ dt \\
&= \T\Big(\Phi^*\big[ t_\Phi\phi(2^k)  \big({\bf 1}-e_{k,N} \big)\big]\Big) + \T\big( \Phi( t_\Phi^{-1}|x_N|) \big)\\
&=\Phi^*\big( t_\Phi\phi(2^k) \big) \T \big({\bf 1}-e_{k,N} \big) + \T\big( \Phi( t_\Phi^{-1}|x_N|) \big)\\
&\leq (2q)^{-1} \Phi^*\big(\phi(2^k) \big) \T \big({\bf 1}-e_{k,N} \big) + \T\big( \Phi( t_\Phi^{-1}|x_N|) \big).
\end{split}
\end{equation*}
where in the last inequality we use \eqref{phi-small}.
According to \cite[p.13]{Kras-Rutickii}, we have 
$u\phi(u)=\Phi(u)+\Phi^*(\phi(u))$. Combining this identity with  \eqref{fact-phi} and  \eqref{e}, we deduce that
\begin{align*}
\Phi(2^k)2^{-k}\T\big(({\bf 1}-e_{k,N})|x_N|\big)&\leq (2q)^{-1} [2^k \varphi(2^k)-\Phi(2^k)] 2^{-k+1} \T \big(({\bf 1}-e_{k,N})|x_N| \big) + \T\big( \Phi( t_\Phi^{-1}|x_N|) \big)\\
&\leq q^{-1} (q-1) \Phi(2^k)2^{-k}\T \big(({\bf 1}-e_{k,N})|x_N| \big) +\T\big( \Phi( t_\Phi^{-1}|x_N|) \big).
\end{align*}
This implies that $\Phi(2^k)2^{-k}\T\big(({\bf 1}-e_{k,N})|x_N|\big)\leq q \T\big( \Phi( t_\Phi^{-1}|x_N|) \big)$. The desired inequality follows from the $\Delta_2$-condition.
\end{proof}
The next lemma may be viewed as  a $\Phi$-moment weak form of   Lemma~\ref{last}.
\begin{lemma}\label{Phi-2}
Let $\Phi$ be an Orlicz  function  that  is $q$-concave for some $1\leq q<2$. Then  there exists a constant $c_{\Phi}'$ so that if $x=(x_n)_{n\geq 1}$ is  a  self-adjoint martingale that is bounded in $L_\Phi(\M)$ and  $k\geq 0$, the following inequality holds for every $N\geq 1$:
\[
\Phi(2^k)2^{-2k}\big\|e_{k,N}x_Ne_{k,N}\big\|_2^2 \leq c_\Phi'  \T\big(\Phi(|x_N|)\big).
\]
\end{lemma}
\begin{proof}
We begin with an estimate recorded  in \eqref{est-1}  that  for $k \geq 0$, we have
\[
\big\|e_{k,N}x_Ne_{k,N}\big\|_2^2 \leq 2 \sum_{j=-\infty}^k 2^{2j}\tau\big({\bf 1}-e_{j-1,N}\big).
\]
Then, from \eqref{e} and Lemma~\ref{Phi-1}, there exists a constant $\a_\Phi$ such that
\[
\Phi(2^k)2^{-2k}\big\|e_{k,N}x_Ne_{k,N}\big\|_2^2 \leq   \a_\Phi \Phi(2^k)2^{-2k}\sum_{j=-\infty}^k  2^{2j} \frac{1}{\Phi(2^{j-1})}  \tau\big(\Phi(|x_N|)\big).
\]
Since  $ t \mapsto t^{-q} \Phi(t)$ is a non-increasing function, it follows that
\[
\Phi(2^k)2^{-2k}\big\|e_{k,N}x_Ne_{k,N}\big\|_2^2 \leq   \a_\Phi 2^{(q-2)k } \sum_{j=-\infty}^k 2^{(2-q)j} \T\big( \Phi(|x_N|)\big).
\]
Since $q<2$, the series is convergent and therefore  the required  estimate is achieved.
\end{proof}

\smallskip

\begin{proof} [Proof of Theorem~\ref{open-B}] As in the proof of Theorem~\ref{main-weak}, we divide the proof into column/row part and diagonal part . We assume that $x$ is a finite martingale $(x_n)_{1\leq n\leq N}$ and consider the same decomposition as in Theorem~\ref{main-weak}. For the column part, we need to verify that for $\lambda>0$,
\begin{equation}\label{Phi-c}
\Phi(\lambda)\T\Big( \ch_{(\lambda,\infty)}\big(\sigma_c(\zeta)\big) \Big) \leq C_\Phi \T\big(\Phi(|x_N|)\big).
\end{equation}
It suffices to verify this for $\lambda=2^k$ for $k\geq 0$. From the estimate on the distribution function stated in  Proposition~\ref{distribution-c}, we have 
\begin{align*}
\Phi(2^k)\T\Big( \ch_{(2^k,\infty)}\big(\sigma_c(\zeta)\big) \Big) &\leq  \Phi(2^k) 2^{-2k} \big\|e_{k,N} x_N e_{k,N} \big\|_2^2  + 8 \Phi(2^k) 2^{-k} \T\big( ({\bf 1}-e_{k,N})|x_N|\big)\\
&\leq \big( c_\Phi' +8 c_\Phi \big) \T\big( \Phi(|x_N|)\big)
\end{align*}
where $c_\Phi$ and $c_\Phi'$ are the constants from Lemma~\ref{Phi-1} and Lemma~\ref{Phi-2} respectively. This proves \eqref{Phi-c}. The row part is identical.

For the diagonal part, we need to verify
\begin{equation}\label{Phi-d}
 \Phi(\lambda)\Tr\big( \ch_{(\lambda, \infty)}(|\eta|) \big) \leq  \T\big( \Phi(|x_N|)\big).
\end{equation}
For $\lambda \geq 1$, it suffices as before to take $\lambda=2^k$ for $k\geq 0$.  This 
can  be deduced as in the column part using Proposition~\ref{distribution-d}(i). Indeed,
\begin{align*}
\Phi(2^k) \Tr\big( \ch_{(2^k, \infty)}(|\eta|) \big)  &\leq
\Phi(2^k)
 2^{-2k+2} \big\|e_{k,N} x_N e_{k,N} \big\|_2^2  + 28 \Phi(2^k) 2^{-k} \T\big( ({\bf 1}-e_{k,N})|x_N| \big)\\
 &+ \Phi(2^k)\T\big(\ch_{(2^k, \infty)}(|x_1|)\big) \\
&\leq \big(4c_\Phi' + 28c_\Phi)\T\big( \Phi(|x_N|)\big) +\Phi(2^k)\T\big(\ch_{(\Phi(2^k), \infty)}(\Phi(|x_1|))\big)\\
&\leq  \big(4c_\Phi' + 28c_\Phi +1)\T\big( \Phi(|x_N|)\big).
\end{align*}
For $0<\lambda<1$, we get from Proposition~\ref{distribution-d}(ii) and  Lemma~\ref{Phi-1} that
\begin{align*}
\Phi(\lambda)\Tr\big( \ch_{(\lambda, \infty)}(|\eta|) \big) &\leq  \Phi(\lambda) \T\big( \ch_{(\lambda,\infty)}(|x_1|) \big) + \Phi(\lambda) \sum_{i\geq 0} \T \big({\bf 1}-e_{i,N} \big)\\
&\leq \Phi(\lambda) \T\big( \ch_{(\Phi(\lambda),\infty)}(\Phi(|x_1|)) \big) + \Phi(\lambda) \sum_{i\geq 0} 2^{-i+1}\T \big( ({\bf 1}-e_{i,N})|x_N| \big)\\
&\leq \T\big(\Phi(|x_N|) \big) + 2c_\Phi \Phi(\lambda)\sum_{i\geq 0} \frac{1}{\Phi(2^i)} \T\big(\Phi(|x_N|) \big).
\end{align*}
Since $t \mapsto t^{-1}\Phi(t)$ is increasing, we deduce that
\[
\Phi(\lambda)\Tr\big( \ch_{(\lambda, \infty)}(|\eta|) \big) \leq  \big( 1 +2c_\Phi \frac{\Phi(\lambda)}{\Phi(1)}\sum_{i\geq 0} 2^{-i}\big)\T\big(\Phi(|x_N|) \big) \leq \big(1 +4c_\Phi\big)  \T\big(\Phi(|x_N|) \big)
\]
where in the second inequality we use the fact that $\Phi$ is an  increasing function.
The proof is complete.
\end{proof}

\begin{proof}[Proof of Theorem~\ref{open-BG}]
In light of Remark~\ref{distribution-similar}, the proof is identical to  the column/row part of Theorem~\ref{open-B}. Details are left to the reader.
\end{proof}

%%%%%%
\subsection{Concluding remarks} 
Assume that $1<p<2$ and $y=(y_n)_{n\geq 1}$ is a  self-adjoint noncommutative martingale  that is weakly differentially subordinate to  another martingale $x=(x_n)_{n\geq 1}$.  From Theorem~\ref{strong-type-1} and Corollary~\ref{H-h}, we have
\begin{equation}\label{h_p}
\big\| y \big\|_{\h_p} \leq c_p \big\|x\big\|_{\h_p} 
\end{equation}
with $c_p =O((p-1)^{-1})$ when $p\to 1$. Similarly, Theorem~\ref{strong-type-2} and Remark~\ref{L-H} yield
\begin{equation}\label{H_p}
\big\| y \big\|_{\H_p} \leq c_p' \big\|x\big\|_{\H_p} 
\end{equation}
with $c_p' =O((p-1)^{-1})$ when $p\to 1$. Inequalities \eqref{h_p} and \eqref{H_p}  are extensions of \eqref{Hardy-1} to the noncommutative setting when $1<p<2$. The fact that the constant in \eqref{Hardy-1} is equal to $1$ naturally leads to the following two questions:
\begin{problem}
Let  $c_p$  and $c_p'$  denote the best constants for   \eqref{h_p} and \eqref{H_p}  respectively. Does there exist an absolute constant  $c$ satisfying 
$c_p \leq c$ and $c_p'\leq c$ for $1<p\leq 2$?
\end{problem}
We note from the Davis decomposition (\cite{Junge-Perrin,Ran-Wu-Xu})  and Corollary~\ref{H-h} that  the set $\{c_p: 1<p\leq 2\}$  is bounded if and only if $\{c_p': 1<p\leq 2\}$ is bounded. 
A related question deals with the case $p=1$.
\begin{problem}\label{H_1} Does there exist a constant $C$ so that  if  $y$ is a self-adjoint martingale that is  weakly differentially subordinate to  another self-adjoint martingale $x\in \H_1(\M)$, then 
\[
\big\| y \big\|_{\H_1} \leq C \big\|x \big\|_{\H_1}?
\]
\end{problem}
Due to the isomorphism between $h_1(\M)$ and $\H_1(\M)$, Problem~\ref{H_1} is equivalent to asking if $\big\| y\big\|_{\h_1} \leq C\big\|x \big\|_{\h_1}$. A positive answer to Problem~\ref{H_1} would imply that if $x \in \H_1(\M)$ then $y$ is an  $L_1$-bounded martingale  thus providing a sufficient condition for the question of when  a weakly differentially subordinate martingale is bounded in $L_1(\M)$. 

\medskip

We conclude the paper with a note on the case $p\geq 2$. Comparing Hardy space norms in relation with differential subordinations  for the case $p\geq 2$ is not as interesting as  the case $1\leq p<2$. This is due to the fact that no decomposition is required for this range and therefore we have the same trivial comparisons of  square functions as in the classical case. Indeed, if  $2\leq p\leq \infty$, $x$ is a self-adjoint  $L_p$-bounded martingale, and $z$ is a self-adjoint martingale that is very weakly differentially subordinate to $x$ in the sense of Definition~\ref{subordinate2}, then  clearly, we have $\|z\|_{\H_p} \leq \|x\|_{\H_p}$ and $\|z\|_{\h_p} \leq \|x\|_{\h_p}$. Moreover, if  $\rm{BMO}(\M)$ is the noncommutatine martingale $\rm{BMO}$-space
(see \cite{PX} for definition) and  
  $x \in \rm{BMO}(\M)$,  then 
$\|z\|_{\rm{BMO}(\M)} \leq  \|x\|_{\rm{BMO}(\M)}$.
%%%%%%

\bigskip 

 \noindent{\bf Acknowledgments.} A portion of the work reported in this paper was carried out  when the second named author visited  Central South University  during the Summer of 2018. It is his pleasure to express his gratitude to all those who made this stay possible  and to the School of Mathematics and Statistics of  the Central South University  for its warm hospitality and financial support.

% \bibliography{narciref,narciref2}
%\bibliographystyle{amsplain}

\def\cprime{$'$}
\providecommand{\bysame}{\leavevmode\hbox to3em{\hrulefill}\thinspace}
\providecommand{\MR}{\relax\ifhmode\unskip\space\fi MR }
\providecommand{\MRhref}[2]{%
  \href{http://www.ams.org/mathscinet-getitem?mr=#1}{#2}
}
\providecommand{\href}[2]{#2}

\end{document}